\documentstyle[12pt,fleqn,leqno]{article}
\title{{\LARGE\bf Partial orderings with the weak Freese-Nation 
property}\\} 
\author{Saka\'e Fuchino$\mbox{}^{\scriptscriptstyle\,0)}$, 
Sabine Koppelberg$\mbox{}^{\scriptscriptstyle\,0)}$ 
and Saharon Shelah$\mbox{}^{\scriptscriptstyle\,1),\,2)}$
\ifcommented
\bigskip\\
\framebox[10.4cm][c]{%
\parbox{10cm}{\footnotesize This is an extended 
version of the paper 
with the same title. Some additional remarks and details of proofs 
are added which will not be in the version for publication.
To remove these remarks, delete the line 
`{\tt $\setminus$commentedtrue \%\%\%}' in the 
LaTeX source file.\smallskip\\
Any comments are appreciated.
}}
\fi
}
\date{{\normalsize November 14, 1994\protect\\(version for submission)%
\protect\protect\\\mbox{}\protect\\
August 24, 1995\protect\\(revised according to some suggestions by the 
referee)}\protect\\\mbox{}\protect\\
\protect\framebox{\protect\framebox{\,to appear in Annals of Pure and Applied Logic\,}}}
\newif\iftesting
\newif\ifcommented 
\let\nc\newcommand
\let\rnc\renewcommand
\let\nth\newtheorem
\let\nenv\newenvironment
\iftesting
\let\Label\label%
\def\label#1{\mbox{}\marginpar{{\tiny #1}}\Label{#1}}%
\fi
\ifcommented
\fi
\nenv{comments}{\ifcommented\begin{footnotesize}}%
{\ifcommented\end{footnotesize}\fi}
\nc{\?}{\iftesting%
\marginpar{$\emptyset\emptyset\emptyset$}%
\fi
\typeout{Something is still incomplete here.}
\errmessage{................................}}
\ifx\kanjiskip\otankonasu
\else\fi
\rnc{\baselinestretch}{1.17}
\nth{Thm}{{\bf Theorem}}[section]
\nth{Cor}[Thm]{{\bf Corollary}}
\nth{Lemma}[Thm]{{\bf Lemma}}
\nth{Prop}[Thm]{{\bf Proposition}}
\nth{Claim}{{\bf Claim}}[Thm]
\nth{Subclaim}{{\bf Subclaim}}[Claim]
\nth{Problem}[Thm]{{\bf Problem}}
\nth{Ex}[Thm]{{\bf Example}}
\nth{Fact}[Thm]{{\bf Fact}}
\nc{\Thmof}[1]{{Theorem \ref{#1}}}
\nc{\Corof}[1]{{Corollary \ref{#1}}}
\nc{\Propof}[1]{{Proposition \ref{#1}}}
\nc{\Lemmaof}[1]{{Lemma \ref{#1}}}
\nc{\Factof}[1]{{Fact \ref{#1}}}
\nc{\Claimof}[1]{{Claim \ref{#1}}}
\nc{\Exof}[1]{{Example \ref{#1}}}
\nc{\Problemof}[1]{{Problem \ref{#1}}}
\nc{\Thmabove}{{Theorem \number\theThm}}
\nc{\Corabove}{{Corollary \number\theThm}}
\nc{\Propabove}{{Proposition \number\theThm}}
\nc{\Lemmaabove}{{Lemma \number\theThm}}
\nc{\Factabove}{{Fact \number\theThm}}
\nc{\Claimabove}{{Claim \number\theClaim}}
\nc{\Exabove}{{Example \number\theThm}}
\nc{\prf}{{\bf Proof\ \ }}
\nc{\prfof}[1]{{\bf Proof of #1\ \ }}
\nc{\prfofClaim}{\raisebox{-.4ex}{\Large $\vdash$\ \ }}
\newsavebox{\qedbox}\sbox{\qedbox}{
{\unitlength=0.034ex \begin{picture}(40,60)
\put(0,0){\framebox(30,44)[cc]{}}
\put(30,-7){\rule{7\unitlength}{44\unitlength}}
\put(10,-7){\rule{27\unitlength}{7\unitlength}}
\end{picture}}}
\nc{\qed}{\nolinebreak\mbox{}\nolinebreak\hfill%
\ifvmode\mbox{}\hfill\fi\usebox{\qedbox}}
\nc{\smallqed}%
{\mbox{}\smallskip\hfill\raisebox{-.4ex}{\Large $\dashv$}\\}
\nc{\qedof}[1]%
{\nopagebreak\mbox{} \hspace*{\fill}{\usebox{\qedbox}{~(#1)}}}
\nc{\Qedof}[1]%
{\nopagebreak\mbox{} \hspace*{\fill}{\usebox{\qedbox}%
{~(#1~\number\theThm)}}}
\nc{\qedofThm}{\Qedof{Theorem}}
\nc{\qedofCor}{\Qedof{Corollary}}
\nc{\qedofProp}{\Qedof{Proposition}}
\nc{\qedofLemma}{\Qedof{Lemma}}
\nc{\qedofFact}{\Qedof{Fact}}
\nc{\qedofClaim}%
{\nopagebreak\mbox{}\hfill\raisebox{-.4ex}{\Large $\dashv$ }\nolinebreak%
\mbox{~(Claim~\number\theClaim)}}
\nc{\qedofSubclaim}%
{\mbox{}\hfill\raisebox{-.4ex}{\Large $\dashv$ }\nolinebreak%
\mbox{~(Subclaim~\number\theSubclaim)}}
\nc{\noindentafterqed}{\bigskip\\}
\nc{\indentafterqed}{\bigskip\\\hspace*{\parindent}}
\nc{\xitem}[1]{\item[$(\mbox{\it #1})$]}
\nc{\xitemof}[1]{{$(\mbox{\it #1})$}}
\nc{\assert}[1]{\makebox[5ex]{${\it #1})$}}
\nc{\lassert}[1]{\llap{\makebox[5ex][l]{${\it #1})$}}}
\nc{\rassert}[1]{\llap{\makebox[5ex][r]{${\it #1})\ \ \ $}}}
\nc{\assertof}[1]{${\it #1})$}
\nenv{assertion}[1]{\begin{trivlist}
\newbox\assertbox
\dimen255=\textwidth
\setbox\assertbox=\hbox{\hspace*{\parindent}{#1}\hspace{\labelsep}}
\advance\dimen255 by -1\wd\assertbox
\advance\dimen255 by -1ex
\item[]\unhbox\assertbox\hfill
\begin{minipage}[t]{\dimen255}}%
{\end{minipage}\end{trivlist}}
\nenv{subassertion}[1]{\begin{trivlist}
\newbox\assertbox
\dimen255=\textwidth
\setbox\assertbox=\hbox{\hspace*{2\parindent}{#1}\hspace{\labelsep}}
\advance\dimen255 by -1\wd\assertbox
\advance\dimen255 by -1ex
\item[]\unhbox\assertbox\hfill
\begin{minipage}[t]{\dimen255}}%
{\end{minipage}\end{trivlist}}
\nc{\implies}{$\,\Rightarrow\,$}
\nc{\equivto}{$\,\Leftrightarrow\,$}
\nenv{formula}[1]{\\[\abovedisplayskip]${#1}$\hfill$}%
{$\hfill\mbox{}\\[\belowdisplayskip]}
\nc{\restr}%
{\mathop{\hspace{0.01ex}|\hspace*{-0.02ex}{\grave{}}\hspace{0.4ex}}}
\nc{\upper}{\mathop{\hspace{0.1ex}\uparrow\hspace{0.2ex}}}
\nc{\subsetnoneq}%
{\mathrel{\raisebox{-0.8ex}{$\stackrel{\subset}{\scriptstyle\,\not=\,}$}}}
\nc{\Linf}[1]{{\cal L}_{\infty{#1}}}
\nc{\powersetof}[1]{{\cal P}(#1)}
\nc{\cardof}[1]{{\mathopen{\mid}{#1}\mathclose{\mid}}}
\nc{\generatedby}[1]{{\langle#1\rangle}}
\nc{\devides}{^{^{\,}}\mid^{^{\,}}}
\nc{\setof}[2]{\{\,#1\,:\,#2\,\}}
\nc{\smallsetof}[1]{\{\,#1\,\}}
\nc{\orderedsetof}[2]{\langle\,#1\,:\,#2\,\rangle}
\nc{\tuple}[1]{\langle\,#1\,\rangle}
\nc{\continuum}{2^{\aleph_0}}
\nc{\reals}{\mbox{$\rm I_{\!\!}R$}}
\nc{\rationals}{%
\rlap{\hbox{\hspace*{0.35ex}\raise0.35ex\hbox{$\scriptstyle\wr$}}}%
\hbox{\rm Q}}
\nc{\integers}{Z\hspace*{-1.12ex}Z}
\nc{\unitint}{\mbox{$\rm I_{\!\!}I$}}
\nc{\lvalue}{\lbrack\!\lbrack}
\nc{\rvalue}{\rbrack\!\rbrack}
\nc{\Bvalueof}[2]{\lvalue\,#1\,\rvalue^{(#2)}}
\nc{\forces}[2]{\,\|\hspace{-.35ex}\mbox{\sf--}_{\,#1\,}%
\mbox{\rm``}\,#2\,\mbox{\rm''}}
\let\Models\models
\def\models#1{\Models\mbox{\rm``}\,#1\,\mbox{\rm''}}
\nc{\xmbox}[1]{ ${\rm #1}$ }
\nc{\st}{such that}
\nc{\wolog}{without loss of generality}
\nc{\Wolog}{Without loss of generality}
\nc{\wrt}{with respect to}
\nc{\Tfae}{The following are equivalent}
\nc{\tfae}{the following are equivalent}
\nc{\Ba}{Boolean algebra}
\nc{\Bas}{Boolean algebras}
\nc{\cBa}{complete Boolean algebra}
\nc{\cBas}{complete Boolean algebras}
\nc{\po}{partial ordering}
\nc{\pos}{partial orderings}
\nc{\Fr}{{\rm Fr}\,}
\nc{\cof}{\mathop{\rm cof}}
\nc{\cf}{\mathop{\rm cf}}
\nc{\pcf}{\mathop{\rm pcf}}
\nc{\Sub}{{\rm Sub}}
\nc{\pr}{{\rm pr}}
\nc{\proj}{{\rm proj}}
\nc{\wproj}{{\rm wproj}}
\nc{\rc}{{\rm rc}}
\nc{\cm}{{cm}}
\nc{\reg}{{\rm reg}}
\nc{\Card}{{\rm Card}}
\nc{\Ord}{{\rm Ord}}
\nc{\Fn}{{\rm Fn}}
\nc{\dom}{{\rm dom}}
\nc{\rng}{{\rm rng}}
\nc{\cohenalg}[1]%
{\mbox{$\,\raisebox{0.05ex}{\small$\wr$}\!\!_{_{\!\!}}\mbox{\rm C}_{#1}$}}
\nc{\cantorsp}[1]{\mbox{}^{#1^{\mbox{}\!}}2}
\nc{\fnsp}[2]{\mbox{}^{#1^{\mbox{}\!}}#2}
\nc{\mapping}[3]{#1:#2\rightarrow #3}
\nc{\ccc}{{\mbox{\rm ccc}}}
\nc{\ZFC}{{\rm ZFC}}
\nc{\CH}{{\rm CH}}
\nc{\GCH}{{\rm GCH}}
\nc{\MA}{{\rm MA}}
\nc{\MM}{{\rm MM}}
\nc{\MAp}{${\rm MA}^+(\sigma${\it-closed\/}$)$}
\nc{\stick}%
{\mbox{{\hspace{0.4ex}$\mbox{\raisebox{1ex}{$\bullet$}}%
\hspace*{-0.94ex}|$\hspace{0.6ex}}}}
\nc{\Subrc}[1]{\Sub_\rc^{\aleph_0}(#1)}
\nc{\dpr}{\overline{\pr}}
\nc{\calB}{{\cal B}}
\nc{\calC}{{\cal C}}
\nc{\calD}{{\cal D}}
\nc{\calF}{{\cal F}}
\nc{\calG}{{\cal G}}
\nc{\calH}{{\cal H}}
\nc{\calN}{{\cal N}}
\nc{\calO}{{\cal O}}
\nc{\calP}{{\cal P}}
\nc{\calS}{{\cal S}}
\nc{\calT}{{\cal T}}
\nc{\dota}{{\dot{a}}}
\nc{\dotb}{{\dot{b}}}
\nc{\dotc}{{\dot{c}}}
\nc{\dotf}{{\dot{f}}}
\nc{\dotg}{{\dot{g}}}
\nc{\dotm}{{\dot{m}}}
\nc{\dotr}{{\dot{r}}}
\nc{\dotx}{{\dot{x}}}
\nc{\doty}{{\dot{y}}}
\nc{\dotB}{{\dot{B}}}
\nc{\dotP}{{\dot{P}}}
\nc{\dotQ}{{\dot{Q}}}
\nc{\dotR}{{\dot{R}}}
\nc{\dotS}{{\dot{S}}}
\nc{\dotU}{{\dot{U}}}
\nc{\dotX}{{\dot{X}}}
\nc{\dotY}{{\dot{Y}}}
\nc{\undercirc}[1]{%
\setbox255=\hbox{$#1$}
\dimen255=0.9ex
\advance\dimen255 by \dp255
\rlap{\hbox to\wd255{\hss{\lower\dimen255\hbox{$\scriptstyle\circ$}}\hss}}
\box255
}
\nc{\sundercirc}[1]{%
\setbox255=\hbox{$\scriptstyle#1$}
\dimen255=0.8ex
\advance\dimen255 by \dp255
\rlap{\hbox to\wd255{\hss{\lower\dimen255
\hbox{$\scriptscriptstyle\circ$}}\hss}}
\box255
}
\nc{\multiundercirc}[1]{%
\mathchoice{\undercirc{#1}}{\undercirc{#1}}{\sundercirc{#1}}{\sundercirc{#1}}}
\nc{\acirc}{\multiundercirc{a}}
\nc{\Dag}{^\dagger}
\nc{\deq}{\mathrel{\sqsubseteq}}
\nc{\Deq}{\mathrel{{\sqsubseteq}\hspace{-1.5ex}{\sqsubseteq}}}
\nc{\notDeq}{\mathrel{{\sqsubseteq}\hspace{-1.5ex}{\not\sqsubseteq}}}
\begin{document}
\maketitle
\rnc{\thefootnote}{\relax}
\footnotetext{\rnc{\baselinestretch}{1.0}\footnotesize
\mbox{}\\[-\medskipamount]
$\mbox{}^{0)}$ 
Institut f\"ur Mathematik II, Freie Universit\"at Berlin, Arnimallee 3, 
14195 Berlin, 
Germany.\\
$\mbox{}^{1)}$ 
Institute of Mathematics, The Hebrew University of Jerusalem, 91904 
Jerusalem, Israel.\\
$\mbox{}^{2)}$ 
Department of Mathematics, Rutgers University, New Brunswick, NJ 08854, U.S.A.

This paper was written while the first author was at the Hebrew University of 
Jerusalem. He would like to thank The Israel Academy of Science and 
Humanities for enabling his stay there. 

The second author was partially supported by the Deutsche 
Forschungsgemeinschaft (DFG) Grant No.\ Ko 490/7--1. 

The third author would like to thank ``Basic Research Foundation'' of The 
Israel Academy of Science and Humanities and the Edmund Landau Center for 
research in Mathematical Analysis supported by the Minerva Foundation 
(Germany) for their support. 
The present paper is the third author's Publication No.\ 549.

The authors would like to thank the referee for some variable suggestions.}
\begin{abstract}
A partial ordering  $P$ is said to have the weak Freese-Nation  
property (WFN) if there 
is a mapping $\mapping{f}{P}{[P]^{\leq\aleph_0}}$ \st, for any $a$, $b\in P$, 
if $a\leq b$ then there exists $c\in f(a)\cap f(b)$ \st\ $a\leq c\leq b$. 
In this note, we study the WFN and some of its generalizations. 

Some features of the class of \Bas\ with the WFN seem to be quite sensitive 
to additional axioms 
of set theory: e.g.,  under \CH, every \ccc\ \cBa\ has this property 
while, 
under ${\bf b}\geq\aleph_2$, there exists no \cBa\ with the WFN. 
(\Thmof{complete-Ba}). 
\end{abstract}
\section{Introduction}\label{introduction}
In \cite{heindorf-shapiro}, a \Ba\ $A$ is said to have the Freese-Nation 
property
(FN, for short) if there exists an FN-mapping 
on $A$, i.e.\ a function $\mapping{f}{A}{[A]^{<\aleph_0}}$ \st\ 
\begin{assertion}{$(*)$}
\it if $a,b\in A$ satisfy $a\leq b$, then $a\leq c\leq b$ for some 
$c\in f(a)\cap f(b)$.
\end{assertion}
This property is closely connected to the notion of freeness because of the 
following 
facts: \assertof{a} {\it every free \Ba\ $A$ has the FN}\/; to see 
this, fix a subset $U$ of $A$ generating $A$ freely; then for $b\in A$, let 
$u(b)$ be a finite subset of $U$ generating $b$ and $f(b)$ the finite 
subalgebra of $A$ generated by $u(b)$. The Interpolation Theorem of 
propositional logic 
then tells us that $(*)$ holds for this $f$. Moreover, 
we have: \assertof{b} 
{\it if $A$ has the FN, then so has every retract of $A$}\/ (see 
\Lemmaof{retract} below for a more general statement). From \assertof{a} and 
\assertof{b}, it follows that: \assertof{c}
{\it every projective \Ba\ has the FN.} 

Historically, the FN was first considered by R.\ Freese and 
J.B.\ Nation in their paper \cite{freese-nation} which gives a 
characterization of projective lattices. In particular, they proved that 
every projective lattice has the FN. 

The FN alone, however, is not equivalent to projectiveness, since, as 
Heindorf proved in \cite{heindorf-shapiro}, a \Ba\ $A$ has the FN 
if and only if $A$ is openly generated in the terminology given below (which 
is also used in \cite{fuchino2}; in \cite{heindorf-shapiro} these 
\Bas\ are called ``rc-filtered''). The notion of open generatedness 
was introduced originally in a topological setting by \v{S}\v{c}epin 
\cite{scepin}. 
In the language of \Bas, a \Ba\ $A$ is said to be openly generated if 
there exists a closed unbounded subset $\calC$ of $[A]^{\aleph_0}$ \st\ 
every $C\in\calC$ is a relatively complete subalgebra of $A$. \v{S}\v{c}epin 
found examples of openly generated \Bas\ which are not projective.

In this paper, we continue the study of \Bas\ with the following weakening of 
the Freese-Nation property, begun in \cite{heindorf-shapiro} or, to 
some extent, already in \cite{scepin}: a \Ba\ $A$ is said to have 
the {\em weak Freese-Nation property} ({\em WFN}\/ for short) if there is a 
{\em WFN mapping on $A$}, that is, a mapping 
$\mapping{f}{A}{[A]^{\leq\aleph_0}}$ satisfying the condition $(*)$ above. 
We solve some open 
problems from \cite{heindorf-shapiro} in Sections \ref{card-functions} and 
\ref{Intalg}. 

Clearly the WFN makes 
perfect sense for arbitrary partial orderings and can be also generalized to 
any uncountable cardinal $\kappa$: we say that a structure $A$ with a 
distinguished partial ordering $\leq$ (we shall call such $A$ a 
{\em partially ordered structure}\/) has the {\em $\kappa$-FN}\/ if there is a 
{\em $\kappa$-FN mapping on $A$}, that is, 
a mapping $\mapping{f}{A}{[A]^{<\kappa}}$ satisfying the condition $(*)$. 
In particular, the FN is the $\aleph_0$-FN and the WFN is the $\aleph_1$-FN.
This generalization is also considered in the following sections.

The paper is organized as follows.  In Section \ref{prelim}, we collect some 
basic facts on the $\kappa$-FN and its connection to the 
$\kappa$-embedding relation $A\leq_\kappa B$ of partially ordered structures.
In Section \ref{chars}, we give some conditions equivalent to the 
$\kappa$-FN which are formulated in terms of elementary submodels, and 
existence of winning strategies in certain infinitary games respectively. 
The behavior of \Bas\ with the $\kappa$-FN \wrt\ the cardinal 
functions of independence, length and cellularity is studied in Section 
\ref{card-functions}. In Sections \ref{Intalg} through 
\ref{Linftykappa-free}, we deal with the question 
which members of the following classes of \Bas\ have 
the WFN: interval algebras, power set algebras, complete \Bas\ 
and $L_{\infty\kappa}$-free \Bas. 

Our notation is standard. For unexplained  
notation and definitions on \Bas, the reader may 
consult \cite{koppelbook} and \cite{monk}.
Some set theoretic notions and basic facts used here can be found in 
\cite{greenbook} and/or \cite{multiple}. 

The authors would like to thank L.\ Heindorf for drawing their attention to 
the weak Freese-Nation property.
\section{$\kappa$-Freese-Nation property and 
$\kappa$-embedding of partially 
ordered structures}\label{prelim}
In this section, we shall look at some basic properties of partially ordered 
structures with 
the $\kappa$-FN. 
In the following, $A$, $B$, $C$ etc.\ are always partially ordered 
structures for 
an arbitrary (but fixed) signature. Note that this setting includes the cases 
that $A$, $B$, $C$ etc.\ are \assertof{a} \Bas\ (with their canonical 
ordering) or \assertof{b} bare partially ordered sets without any additional 
structure. 

By the theorem of Heindorf mentioned above, 
every openly generated \Ba\ has the WFN. But the class of \Bas\ 
with the WFN contains many  
more \Bas. This can be seen already in the following:
\begin{Lemma}\label{aleph-1}
If $\cardof{A}\leq\kappa$ then $A$ has the $\kappa$-FN. 
\end{Lemma}
\prf Let $A=\setof{b_\alpha}{\alpha<\kappa}$. The mapping 
$\mapping{f}{A}{[A]^{<\kappa}}$ defined by 
\[ f(b_\alpha)=\setof{b_\beta}{\beta\leq\alpha}
\]\noindent
for $\alpha<\kappa$, is a $\kappa$-FN mapping on $A$.\qedofLemma
\indentafterqed
For $A$, $B$ \st\ $A\leq B$ (i.e.\/ $A$ is a substructure of 
$B$) and $b\in B$ we write:
\[ 
\begin{array}{l}
A\restr b =\setof{a\in A}{a\leq b},\medskip\\
A\upper b =\setof{a\in A}{a\geq b}. 
\end{array}
\]\noindent
$A$ is a {\em$\kappa$-substructure of $B$} (or {\em$\kappa$-subalgebra of 
$B$} in case of \Bas; notation: $A\leq_\kappa B$) if 
$A\leq B$ and, for every $b\in B$, there are a cofinal subset $U$ of 
$A\restr b$ 
and a coinitial subset $V$ of $A\upper b$ both of cardinality less than 
$\kappa$. 
For $\kappa=\aleph_1$ we say also that $A$ is a 
{\em$\sigma$-substructure/subalgebra} 
of $B$ and denote it by $A\leq_\sigma B$. For $\kappa=\aleph_0$, 
a $\kappa$-substructure/subalgebra $A$ of $B$ is also called a {\em relatively 
complete substructure/subalgebra}\/) of $B$ and this is denoted also by 
$A\leq_\rc B$. Note that, if 
$\leq$ is 
lattice order on $A$, then $A\leq_\rc B$ holds if and only if, 
for all $b\in B$, 
$A\restr b$ has a cofinal subset $U$ and $A\upper b$ has a coinitial subset 
$V$ consisting of a single element respectively. In this case these elements 
are called the {\em lower} and the {\em upper projection of $b$ on $A$} and 
denoted by 
$p^B_A(b)$ and $q^B_A(b)$ respectively. 
Note also that, for \Bas, to show that $A\leq_\kappa B$ holds, it is 
enough to check that $A\restr b$ is $<\kappa$-generated for every $b\in B$, 
by duality.

The following lemma can be proved easily:
\begin{Lemma}\label{trans}
\assert{a} If $\lambda\leq\kappa$ and $A\leq_\lambda B$ then 
$A\leq_\kappa B$. In particular,  
if $A\leq_\rc B$ then $A\leq_\sigma B$.\\
\assert{b} If $A\leq_\kappa C$ and $A\leq B\leq C$ then 
$A\leq_\kappa B$.\\
\assert{c}
For a regular cardinal $\kappa$, if $A\leq_\kappa B$ and $B\leq_\kappa C$ then 
$A\leq_\kappa C$.\qed 
\end{Lemma}
\begin{Lemma}\label{sigma-embed}
\assert{a} For a regular cardinal $\kappa$, if $B$ has the $\kappa$-FN 
and $A\leq_\kappa B$ then $A$ also 
has the $\kappa$-FN.\smallskip\\
\assert{b} For a regular cardinal $\kappa$, if $A\leq_\kappa B$, $B$ has the 
$\kappa$-FN and $f$ is a $\kappa$-FN 
mapping on $A$, then there is a $\kappa$-FN mapping $\tilde{g}$ on $B$ 
extending $f$. \smallskip\\
\assert{c} If ${g}$ is a $\kappa$-FN mapping on $B$ and $C\leq B$ is 
closed \wrt\ $g$ $($i.e.\ $g(c)\subseteq C$ holds for all $c\in C$\/$)$, then 
$C\leq_\kappa B$. 
\end{Lemma} 
\prf
For \assertof{a} and \assertof{b}, 
let $\mapping{g}{B}{[B]^{<\kappa}}$ be a $\kappa$-FN mapping on $B$ and,  
for each $b\in B$, let 
$U(b)$ and $V(b)$ be \st\ 
$U(b)$ is a cofinal subset of $A\restr b$, $V(b)$ is a coinitial subset 
of $A\upper b$ and $\cardof{U(b)},\cardof{V(b)}<\kappa$.\smallskip\\
\assertof{a}: 
Let $f$ be the mapping on $A$ defined by
\[ f(a)=\bigcup\setof{U(b)}{b\in g(a)}.
\]\noindent
Since $\kappa$ is regular we have $f(a)\in[A]^{<\kappa}$ for every $a\in A$. 
$f$ is a $\kappa$-FN mapping on $A$.
To see this let $a,a'\in A$ be \st\ $a\leq a'$. Then there is 
$b\in g(a)\cap g(a')$ \st\ $a\leq b\leq a'$. Since $U(b)$ is cofinal 
in $A\restr b$, there is $c\in U(b)$ ($\subseteq f(a)\cap f(a')$) \st\ 
$a\leq c$. 
Since $c\leq b$ we also have $c\leq a'$. Note that in this proof we 
only needed that one of $U(b)$ and $V(b)$ is of cardinality less than 
$\kappa$ for every $b\in B$.
\smallskip\\
\assertof{b}: Let $\tilde g$ be the mapping on $B$ defined by
\[ \tilde{g}(b)=
\left\{\,
\begin{array}{@{}l@{}l}
f(b) &;\mbox{ if }b\in A,\\
g(b)\cup\bigcup\setof{f(c)}{c\in U(b)\cup V(b)}\qquad&; 
\mbox{ otherwise.}
\end{array}
\right.
\]\noindent
Clearly $f\subseteq\tilde{g}$. 
$\tilde{g}$ is a $\kappa$-FN mapping: since $\kappa$ is regular, we have
$\tilde{g}(b)\in[B]^{<\kappa}$ for every $b\in B$. Let $b,b'\in B$ be 
\st\ $b\leq b'$. We 
want to show that there is $c\in \tilde{g}(b)\cap\tilde{g}(b')$ \st\ 
$b\leq c\leq b'$. If $b$, $b'\in A$ or $b$, $b'\in B\setminus A$, this 
follows immediately from the definition of $\tilde{g}$. Suppose that 
$b\in A$ and $b'\in B\setminus A$. Then there is $d\in U(b')$ \st\ 
$b\leq d$. Hence there is $c\in f(b)\cap f(d)$ \st\ $b\leq c\leq d$ holds. 
Since $f(b)=\tilde{g}(b)$ and $f(d)\subseteq g(b')$ by $d\in U(b')$, it 
follows that $c\in \tilde{g}(b)\cap\tilde{g}(b')$ and $b\leq c\leq b'$. The 
case, $b\in B\setminus A$ and $b'\in A$, can be treated similarly. 
\smallskip\\
\assertof{c}: Let $C\leq B$ be closed \wrt\ $g$.
For $b\in B$, let $U=g(b)\cap(C\restr b)$ and $V=g(b)\cap(C\upper b)$. Then 
clearly $U$ and $V$ are of cardinality $<\kappa$. 
We show that $U$ is cofinal in $C\restr b$\,: 
if $c\leq b$ for some $c\in C$ then there is 
$e\in g(c)\cap g(b)$ \st\ $c\leq e\leq b$ holds. Since $g(c)\subseteq C$, 
we have $e\in C\restr b$. Hence $e\in U$. Similarly we can also show that 
$V$ is coinitial in $C\upper b$.
\qedofLemma\smallskip\\

As already mentioned in the introduction, the \Bas\ with the FN property are 
exactly the openly generated \Bas\ (\cite{heindorf-shapiro}). Hence it 
follows from the next lemma that, if $(B_\alpha)_{\alpha<\delta}$ is a 
continuously increasing chain of openly generated \Bas\ \st\ 
$B_\alpha\leq_\sigma B_{\alpha+1}$ for every $\alpha<\delta$, then 
$\bigcup_{\alpha<\delta}B_\alpha$ is also openly generated. 
The original proof of this fact in \cite{scepin} employed very 
complicated combinatorial arguments, while our proof below and also the 
proof of the characterization of openly generated \Bas\ as those with the FN 
property is quite elementary. 
\begin{Lemma}\label{conti-chain}Suppose that $\kappa$ is a regular cardinal,
$\delta$ a limit ordinal and $(B_\alpha)_{\alpha\leq\delta}$ a 
continuously increasing chain \st\ 
$B_\alpha\leq_\kappa B_{\alpha+1}$ for all $\alpha<\delta$. \smallskip Then 
\ifvmode\\\fi
\assert{a} $B_\alpha\leq_\kappa B_\beta$ for every 
$\alpha\leq\beta\leq\delta$. \smallskip\\
\assert{b}If $B_\alpha$ has the $\kappa$-FN  
for every $\alpha<\delta$, then 
$B_\delta$ also has the $\kappa$-FN.
\end{Lemma}
\prf \assertof{a}: By induction on $\beta$, using 
\Lemmaof{trans},\,\assertof{c} for successor steps.\smallskip\\ 
\assertof{b}: By \Lemmaof{sigma-embed},\,\assertof{b}, we can construct a 
continuously increasing sequence $(f_\alpha)_{\alpha<\delta}$ \st\ 
for each $\alpha<\delta$, $f_\alpha$ is a $\kappa$-FN mapping on $B_\alpha$. 
$f_\delta=\bigcup_{\alpha<\delta}f_\alpha$ is then a $\kappa$-FN mapping on 
$B_\delta$. 
\qedofLemma
\begin{Lemma}\label{countable-union}
Suppose that $\mu<\kappa$, $\cf(\mu)<\cf(\kappa)$ and  
$(B_\alpha)_{\alpha\in\mu}$ is an increasing sequence of 
$\kappa$-substructures of $B$. Then 
$\bigcup_{\alpha\in\mu}B_\alpha$ is also a $\kappa$-substructure of $B$. 
\end{Lemma}
\prf
\Wolog\ we may assume that $\mu=\cf(\mu)$ holds. 
For $b\in B$ let $U_\alpha(b)$ be a cofinal subset of $B_\alpha\restr b$ and 
$V_\alpha(b)$ a coinitial subset of $B\upper b$ both of cardinality less 
than $\kappa$. Then 
$U(b)=\bigcup_{\alpha\in\mu}U_\alpha(b)$ is a cofinal subset of  
$(\bigcup_{\alpha\in\mu}B_\alpha)\restr b$ and 
$V(b)=\bigcup_{\alpha\in\mu}V_\alpha(b)$ is a coinitial subset of 
$(\bigcup_{\alpha\in\mu}B_\alpha)\upper b$. Since $\mu<\cf(\kappa)$, we have 
$\cardof{U(b)}, \cardof{V(b)}<\kappa$.
\qedofLemma
\begin{Lemma}\label{non-conti-chain}
Suppose that $\kappa$ is a regular cardinal, 
$\delta$ a limit ordinal and $(A_\alpha)_{\alpha<\delta}$ an 
increasing chain \st\ 
$A_\alpha\leq_\kappa A_{\alpha+1}$ for all 
$\alpha<\delta$ and $A_\gamma=\bigcup_{\alpha<\gamma}A_\alpha$ for 
all limit $\gamma<\delta$ with $\cf(\gamma)\geq\kappa$. 
Let $A=\bigcup_{\alpha<\delta}A_\alpha$. 
If $A_\alpha$ has the $\kappa$-FN for every $\alpha<\delta$, then 
$A$ also has the $\kappa$-FN.
\end{Lemma}
\prf
Let $(B_\alpha)_{\alpha\leq\delta}$ be defined by:
\[ B_\alpha=
\left\{\,
\begin{array}{@{}ll}
	A_\alpha &;\mbox{ if }\alpha\mbox{ is a successor or of cofinality }
	\geq\kappa,\\
	{\bigcup}_{\beta<\alpha}A_\beta &;\mbox{ otherwise.} 
\end{array}\right.
\]\noindent
Then $(B_\alpha)_{\alpha\leq\delta}$ is continuously increasing, 
$B_\delta=A$ and $B_\alpha\leq_{\kappa}B_{\alpha+1}$ for all $\alpha<\delta$: 
for a limit $\alpha<\delta$ with $\cf(\alpha)<\kappa$, this follows 
from \Lemmaof{countable-union}. 

Hence, using \Lemmaof{conti-chain},\,\assertof{b}, we can show by induction 
that $B_\alpha$ has the $\kappa$-FN for every $\alpha\leq\delta$. 
\qedofLemma
\begin{Lemma}\label{retract}
Suppose that there are order preserving mappings $\mapping{i}{A}{B}$ and 
$\mapping{j}{B}{A}$ \st\ $j\circ i=id_A$. If $B$ has the $\kappa$-FN, then 
$A$ also has the $\kappa$-FN. In particular, for \Bas\ $A$, $B$, if 
$A$ is a retract of $B$ and $B$ has the $\kappa$-FN, then $A$ also has the 
$\kappa$-FN.
\end{Lemma}
\prf
Let $\mapping{g}{B}{[B]^{<\kappa}}$ be a $\kappa$-FN mapping on $B$ and 
$f$ be the mapping on $A$ defined by
\[ f(a)=j[g(i(a))].
\]\noindent
We show that $f$ is a $\kappa$-FN mapping on $A$. Clearly 
$f(a)\in[A]^{<\kappa}$ for every $a\in A$. Suppose that $a,a'\in A$ are  
\st\ $a\leq a'$. Then we have $i(a)\leq i(a')$. Hence there is  
$b\in g(i(a))\cap g(i(a'))$ \st\ $i(a)\leq b\leq i(a')$. It follows that 
\[ a=j\circ i(a)\leq j(b)\leq j\circ i(a')=a'
\]\noindent
and $j(b)\in f(a)\cap f(a')$.
\qedofLemma 
\section{Characterizations of partially ordered structures with the weak 
Freese-Nation property}\label{chars} 
For a partially ordered structure $B$, let us say that a regular cardinal 
$\chi$ is sufficiently large if 
the $n$'th power of $B$ ${\cal P}^n(B)$ is in $\calH_\chi$ for every 
$n\in\omega$, where $\calH_\chi$ is the set of every sets of hereditary of 
cardinality less than $\chi$. 
\begin{Prop}\label{characterization}
For a regular $\kappa$ and a partially ordered structure $B$, 
\tfae:\smallskip\\
\assert{1} $B$ has the $\kappa$-FN;\smallskip\\
\assert{2} For some, or equivalently, any sufficiently large $\chi$, if 
$M\prec\calH_\chi=(\calH_\chi,\in)$ is \st\ 
$B\in M$, $\kappa\subseteq M$ and $\cardof{M}=\kappa$ then 
$B\cap M\leq_\kappa B$ 
holds;\smallskip\\
\assert{3} $\setof{C\in[B]^{\kappa}}{C\leq_\kappa B}$ contains a club set;
\smallskip\\
\assert{4} There exists a \po\ $I=(I,\leq)$ and an indexed family 
$(B_i)_{i\in I}$ of substructures of $B$ of cardinality $\kappa$ \st
\begin{subassertion}{}
\mbox{}\rassert{i}%
$\setof{B_i}{i\in I}$ is cofinal in $([B]^{\kappa},\subseteq)$,\\
\rassert{i{}i}%
$I$ is directed and for any $i$, $j\in I$, if $i\leq j$ then $B_i\leq B_j$,\\
\rassert{i{}i{}i}%
for every well-ordered $I'\subseteq I$ of cofinality $\leq\kappa$, 
$i'=\sup I'$ 
exists and $B_{i'}=\bigcup_{i\in I'}B_i$ holds, and\\
\rassert{i{}v}%
$B_i\leq_\kappa B$ holds for every $i\in I$.
\end{subassertion}
\end{Prop}
\prf
\assertof{1}\implies\assertof{2}: 
Let $f$ be a $\kappa$-FN mapping on $B$. Since $\chi$ is sufficiently large 
for $B$, we 
have $B$, $f\in\calH_\chi$. Let $M\prec H_\chi$ be \st\ $B\in M$, 
$\kappa\subseteq M$ and 
$\cardof{M}=\kappa$. Then
there is a $\kappa$-FN mapping $f'$ on $B$ in $M$. Clearly $B\cap M$ is 
closed \wrt\ 
$f'$. Hence it follows by \Lemmaof{sigma-embed} that 
$B\cap M\leq_\kappa B$. \smallskip\\ 
\assertof{2}\implies\assertof{3}: 
Clear.\smallskip\\
\assertof{3}\implies\assertof{4}: 
Let $I\subseteq\setof{C\in[B]^{\kappa}}{C\leq_\kappa B}$ be a club 
subset of $[B]^{\kappa}$ with the substructure relation. For $A\in I$, let
$B_A=A$. 
Then $(I,\leq)$ and 
$(B_A)_{A\in I}$ satisfy the 
conditions in \assertof{4}.\smallskip\\ 
\assertof{4}\implies\assertof{1}: we prove this in the following two claims.
Let $I$ and $(B_i)_{i\in I}$ be as in \assertof{4}.
For a directed $I'\subseteq I$, let $B_{I'}=\bigcup_{i\in I'}B_i$. 
\begin{Claim}\label{cl0}
If $I'\subseteq I$ is directed, then $B_{I'}\leq_\kappa B$.
\end{Claim}
\prfofClaim
Otherwise there is $b\in B$ \st\ either $B_{I'}\restr b$ does not 
have any cofinal subset of cardinality less than $\kappa$ or 
$B_{I'}\upper b$ does 
not have any coinitial subset of cardinality less than $\kappa$. For 
simplicity, let 
us assume the first case. 
Then there exists an increasing sequence $(I_\alpha)_{\alpha<\kappa}$ of 
directed subsets of $I'$ of cardinality less than $\kappa$ \st\ 
$B_{I_\alpha}\restr b$ is not cofinal in $B_{I_{\alpha+1}}\restr b$. 
By \assertof{iii}, $i_\alpha=\sup I_\alpha$ exists and 
$B_{i_\alpha}=B_{I_\alpha}$ holds for every $\alpha<\kappa$. 
$(i_\alpha)_{\alpha<\kappa}$ is an increasing sequence in $I$. Hence, 
again by \assertof{iii}, there exists $i^*=\sup_{\alpha<\kappa}i_\alpha$ 
and $B_{i^*}=\bigcup_{\alpha<\kappa}B_{i_\alpha}$. By \assertof{iv}, 
$B_{i^*}\leq_\kappa B$. But by the construction, $B_{i^*}\restr b$ cannot 
have any cofinal subset of cardinality less than $\kappa$. This is a 
contradiction. 
\qedofClaim
\begin{Claim}\label{cl1}
If $I'\subseteq I$ is directed, then $B_{I'}$ has the $\kappa$-FN.
\end{Claim}
\prfofClaim
We prove the claim by induction on $\cardof{I'}$. 
If $\cardof{I'}\leq\kappa$, we have $\cardof{B_{I'}}=\kappa$. Hence, by 
\Lemmaof{aleph-1}, $B_{I'}$ has the $\kappa$-FN.
Assume that $\cardof{I'}=\lambda>\kappa$ and that we have 
proved the claim for every directed
$I''\subseteq I$ with $\cardof{I''}<\lambda$.
Take a continuously increasing sequence $(I_\alpha)_{\alpha<\cf(\lambda)}$ of 
directed subsets of $I'$ \st\ $\cardof{I_\alpha}<\lambda$ for every 
$\alpha<\cf(\lambda)$ and $I'=\bigcup_{\alpha<\cf(\lambda)}I_\alpha$. 
$(B_{I_\alpha})_{\alpha<\cf(\lambda)}$ is then a continuously increasing 
sequence 
of substructures of $B_{I'}$ and 
$B_{I'}=\bigcup_{\alpha<\cf(\lambda)}B_{I_\alpha}$. 
By the induction hypothesis, $B_{I_\alpha}$ has the $\kappa$-FN and, 
by 
\Claimof{cl0}, we have $B_{I_\alpha}\leq_\kappa B$. Hence, by 
\Lemmaof{conti-chain},\assertof{b}, $B_{I'}$ has also the $\kappa$-FN. 
\qedofClaim
\noindentafterqed
Now by applying \Claimabove\ to $I'=I$, we can conclude that 
$B=B_I$ has the $\kappa$-FN. 
\qedofProp%
\nc{\kappagame}[2]{\calG^{#1}(#2)}%
\nc{\kappalambdagame}[3]{\calG^{#1}_{#2}(#3)}%
\indentafterqed
Now we give yet another characterization of partially ordered structures 
with the $\kappa$-FN 
by means of a game. This characterization will be used later in 
the proof of Propositions \ref{omega-2}, \ref{ccc-rc}, etc. 
For a partially ordered structure\ $B$, let $\kappagame{\kappa}{B}$ be the 
following game played by Players I and II: 
in a play in $\kappagame{\kappa}{B}$, Players I and II choose 
subsets 
$X_\alpha$ and $Y_\alpha$ of $B$ of cardinality less than $\kappa$ 
alternately for 
$\alpha<\kappa$ \st
\[ 
X_0\subseteq Y_0\subseteq X_1\subseteq Y_1\subseteq\cdots\subseteq X_\alpha
\subseteq Y_\alpha\subseteq\cdots\subseteq X_\beta\subseteq 
Y_\beta\subseteq\cdots 
\]\noindent
for $\alpha\leq\beta<\kappa$. So a play in $\kappagame{\kappa}{B}$ 
looks 
like 
\[ 
\begin{array}{l@{}l@{}@{}l@{}@{}l@{}l}
\mbox{\it Player I}\ &:\ \ &X_0,\ &X_1,\ \ldots,\ &X_\alpha,\ \ldots\medskip\\ 
\mbox{\it Player II}\ &:\ \ &Y_0,\ &Y_1,\ \ldots,\ &Y_\alpha,\ \ldots
\end{array}
\]\noindent
where $\alpha<\kappa$. Player II wins the play if 
$\bigcup_{\alpha<\kappa}X_\alpha=\bigcup_{\alpha<\kappa}Y_\alpha$ is a 
$\kappa$-substructure of $B$. 
Let us call a strategy $\tau$ for Player II simple if, in $\tau$, each 
$Y_\alpha$ is decided from the information of the set 
$X_\alpha\subseteq B$ alone (i.e.\ also independent of $\alpha$). 

For a sufficiently large $\chi$ (\wrt\ $B$), an elementary submodel $M$ of 
$\calH_\chi$ is said to be $V_{\kappa}$-like if, either 
$\kappa=\aleph_0$ and $M$ is countable, or there is an 
increasing 
sequence $(M_\alpha)_{\alpha<\kappa}$ of elementary submodels 
of $M$ of cardinality less than $\kappa$ \st\ 
$M_\alpha\in M_{\alpha+1}$ 
for all 
$\alpha<\kappa$ and $M=\bigcup_{\alpha<\kappa}M_\alpha$. If $M$ is 
$V_{\kappa}$-like, we say that a sequence
$(M_\alpha)_{\alpha<\kappa}$ as above witnesses the $V_{\kappa}$-likeness of 
$M$.
The notion of $V_\kappa$-like elementary submodels of $\calH_\chi$ is a 
weakening of internally approachable elementary submodels introduced in 
\cite{FMSh:240}. 
An elementary submodel $M$ of $\calH_\chi$ is said to be internally 
approachable if $M$ is the union of continuously increasing sequence 
$(M_\alpha)_{\alpha<\kappa}$ of smaller elementary submodels \st\ 
$(M_\beta)_{\beta\leq\alpha}\in M_{\alpha+1}$ for every $\alpha<\kappa$. 
The main reason of the use of $V_\kappa$-like elementary submodels here 
instead of internally approachable ones is the following 
\Lemmaof{V-kappa-like},\,\assertof{b} which seems to be false in general for 
internally approachable elementary submodels.
\begin{Lemma}\label{V-kappa-like}
\assert{a}
If $M$ is a $V_{\kappa}$-like elementary submodel of $\calH_\chi$ \st\ 
$\kappa\in M$, then 
$\kappa\subseteq M$ holds. Hence, if $x$ is of cardinality less or equal to 
$\kappa$ and 
$x\in M$ then we have $x\subseteq M$. 
\medskip\\
\assert{b}
If $(N_\alpha)_{\alpha<\kappa}$ is an increasing sequence of 
$V_{\kappa}$-like elementary submodels of $\calH_\chi$, then 
$M=\bigcup_{\alpha<\kappa}N_\alpha$ is also a $V_{\kappa}$-like elementary 
submodel of $\calH_\chi$.
\end{Lemma}
\prf
\assertof{a}:
Let $(M_\alpha)_{\alpha<\kappa}$ witness the $V_{\kappa}$-likeness of 
$M$. 
Assume that $\kappa\not\subseteq M$. Let
\[ \alpha_0=\min\setof{\alpha\in\kappa}{\alpha\not\in M}.
\]\noindent
Then we have $\alpha_0\subseteq M$. Let
\[ \alpha_1=\min\setof{\alpha\leq\kappa}{\alpha_0\leq\alpha,\,\alpha\in M}.
\]\noindent
Since $\alpha_0$ is of cardinality less than $\kappa$, there exists 
$\alpha<\kappa$ 
\st\ 
$\alpha_0\subseteq M_\alpha$. Then we have $\alpha_1\in M_{\alpha+1}$. 
Since $M_\alpha\in M_{\alpha+1}$, 
$\alpha_0=\setof{\beta\in M_\alpha}{\beta<\alpha_1}$ is an element of 
$M_{\alpha+1}\subseteq M$. This is a contradiction. Hence we have
$\kappa\subseteq M$.

If $x$ is of cardinality less or equal to $\kappa$ and $x\in M$, then there 
is a 
surjection $\mapping{f}{\kappa}{x}$ in $M$. Since $\kappa\subseteq M$, 
it follows that $x=f[\kappa]\subseteq M$.
\medskip\\
\assertof{b}:
It is clear that $M$ is an elementary submodel of $\calH_\chi$. To prove 
that $M$ is $V_{\kappa}$-like, let
$M=\setof{m_\xi}{\xi<\kappa}$ and, 
for each $\alpha<\kappa$, let $(N_{\alpha,\beta})_{\beta<\kappa}$ be an 
increasing sequence of elementary submodels of $\calH_\chi$ of cardinality 
less than $\kappa$ 
witnessing the $V_{\kappa}$-likeness of $N_\alpha$. Since 
$N_{\alpha,\beta}\in N_\alpha\subseteq M$ and 
$\bigcup_{\beta<\kappa}N_{\alpha,\beta}=N_\alpha$, 
we can choose 
$\alpha_\xi$, $\beta_\xi<\kappa$ for $\xi<\kappa$ inductively 
\st\medskip\\ 
\assert{a} $(N_{\alpha_\xi,\beta_\xi})_{\xi<\kappa}$ is an increasing 
sequence,\smallskip\\ 
\assert{b} $N_{\alpha_\xi,\beta_\xi}\in N_{\alpha_{\xi+1},\beta_{\xi+1}}$ 
holds for every $\xi<\kappa$ and\smallskip\\
\assert{c} $m_\xi\in N_{\alpha_{\xi},\beta_{\xi}}$ for every
$\xi<\kappa$.\medskip\\ 
Then $(N_{\alpha_\xi,\beta_\xi})_{\xi<\kappa}$ witnesses the 
$V_{\kappa}$-likeness of $M$.
\qedofLemma 
\begin{Prop}\label{game}
For regular $\kappa$ and a partially ordered structure\ $B$, 
\tfae:\smallskip\\ 
\assert{1} $B$ has the $\kappa$-FN;\smallskip\\
\assert{2} Player II has a simple winning strategy in 
$\kappagame{\kappa}{B}$;\smallskip\\
\assert{3} For some, or equivalently any, sufficiently large $\chi$, if 
$M\prec\calH_\chi$ is $V_{\kappa}$-like \st\ $B$, $\kappa\in M$, 
then $B\cap M\leq_\kappa B$.
\end{Prop}
\prf
Assume that $\kappa$ is uncountable
(for $\kappa=\aleph_0$, the proof is easier than the following one and given 
in \cite{appendix}).\medskip\\
\assertof{1}\implies\assertof{2}: Let $\mapping{f}{B}{[B]^{<\kappa}}$ be a 
$\kappa$-FN mapping on $B$. Then Player II can win by the following strategy: 
in the 
$\alpha$'th move, Player II chooses $Y_\alpha$ so that 
$X_\alpha\subseteq Y_\alpha$ and $Y_\alpha$ is a substructure of 
$B$ of cardinality less than $\kappa$ closed under $f$. After $\kappa$ moves, 
$\bigcup_{\alpha<\kappa}Y_\alpha$ is a substructure of $B$ closed under 
$f$. Hence, by \Lemmaof{sigma-embed},\,\assertof{c}, it is a 
$\kappa$-substructure of $B$.\medskip\\
\assertof{2}\implies\assertof{3}: Let $M$ be a $V_{\kappa}$-like 
elementary submodel of $\calH_\chi$ \st\ $B$, $\kappa\in M$. We have to 
show that 
$B\cap M\leq_\kappa B$. 
Let $(M_\alpha)_{\alpha<\kappa}$ witness the $V_{\kappa}$-likeness of 
$M$. \Wolog\ we may assume that $B\in M_0$. By $M_0\prec\calH_\chi$, there 
is a simple winning strategy 
$\tau\in M_0$ for Player II in $\kappagame{\kappa}{B}$ (hence 
$\tau\in M_\alpha$ for every $\alpha<\kappa$).  Let 
$(X_\alpha, Y_\alpha)_{\alpha<\kappa}$ be the play in 
$\kappagame{\kappa}{B}$ \st\ at his $\alpha$'th move, Player I took 
$B\cap M_{\xi_\alpha}$ for some $\xi_\alpha<\kappa$ and 
Player II played always according to $\tau$. Such a game is possible since 
if Player I chooses $B\cap M_{\xi_\alpha}$ at his $\alpha$'th move, then 
$B\cap M_{\xi_\alpha}\in M_{\xi_\alpha+1}$. Hence Player II's  move $Y_\alpha$ 
taken according to $\tau$ is also an element of $M_{\xi_\alpha+1}$. Since 
$Y_\alpha$ 
is of cardinality less than $\kappa$, we have $Y_\alpha\subseteq M$ by 
\Lemmaof{V-kappa-like},\,\assertof{a}. Hence there is some 
$\xi_{\alpha+1}\geq\xi_\alpha$ \st\ 
$Y_\alpha\subseteq M_{\xi_{\alpha+1}}$. 
Thus Player I may take $B\cap M_{\xi_{\alpha+1}}$ at his next move. 

Now we have 
$B\cap M=B\cap(\bigcup_{\alpha<\kappa}M_{\xi_\alpha})\leq_\kappa B$ 
since $\tau$ was a winning strategy of Player II.\medskip\\
\assertof{3}\implies\assertof{1}: 
Let 
\[ 
\begin{array}{@{}ll}
\calC=\setof{M}{&\cardof{M}=\kappa,\,M\mbox{ is a union of an increasing 
sequence}\\
&\mbox{of }V_{\kappa}\mbox{-like elementary submodels of }
\calH_\chi\mbox{ \st\ }B,\,\kappa\in M}. 
\end{array}
\]\noindent
Then it is easy to see that $\calC$ is club in $[\calH_\chi]^{\kappa}$. 
Hence $\calC'=\setof{B\cap M}{M\in\calC}$ contains a club subset of 
$[B]^{\kappa}$. By \Propof{characterization},\,\assertof{3}, if follows 
from the claim below that $B$ has the $\kappa$-FN.
\begin{Claim}
$B\cap M\leq_\kappa B$ holds for every $M\in \calC$.
\end{Claim}
\prfofClaim
If $M\in \calC$ is union of a sequence $(M_\alpha)_{\alpha<\rho}$ of 
$V_{\kappa}$-like elementary submodels of $\calH_\chi$ for some 
$\rho<\kappa$ \st\ $B$, $\kappa\in M_0$, then we 
have 
$B\cap M_\alpha\leq_\kappa B$ for every $\alpha<\rho$ by \assertof{3}. Hence 
we have $B\cap M=\bigcup_{\alpha<\rho}(B\cap M_\alpha)\leq_\kappa B$ by 
\Lemmaof{countable-union}. If $M\in \calC$ is the union 
of a $\kappa$-chain of $V_{\kappa}$-like elementary submodels of 
$\calH_\chi$, then if follows by \Lemmaof{V-kappa-like},\,\assertof{b} 
that $M$ itself is 
$V_{\kappa}$-like. Hence we have $B\cap M\leq_\kappa B$ by the assumption. 
\qedofClaim
\qedofProp
\indentafterqed
Under $2^{<\kappa}=\kappa$, Propositions \ref{characterization} and 
\ref{game} can be 
yet improved. This is because of the following fact:
\begin{Lemma}\label{Vomega-1-like}
Assume that $\kappa$ is a regular cardinal \st\ $2^{<\kappa}=\kappa$. 
Let $B$ be a partially ordered structure, $\chi$ be sufficiently 
large for $B$ and $M\subseteq\calH_\chi$. 
Then \tfae:\smallskip\\ 
\assert{1} $M$ is a $V_{\kappa}$-like elementary submodel of
$\calH_\chi$;\smallskip\\
\assert{2} $M\prec\calH_\chi$, $\cardof{M}=\kappa$, $B\in M$ and
$[M]^{<\kappa}\subseteq M$.
\end{Lemma}
\prf
For $\kappa=\aleph_0$, this is clear. Assume that $\kappa$ is uncountable. 
\smallskip\\
\assertof{1}\implies\assertof{2}: Let $M$ be a $V_{\kappa}$-like 
elementary submodel of $\calH_\chi$ and let
$(M_\alpha)_{\alpha<\kappa}$ be an increasing sequence of 
elementary submodels of $\calH_\chi$ witnessing the 
$V_{\kappa}$-likeness of $M$. It is enough to show that
$[M_\alpha]^{<\kappa}\subseteq M$ holds for every $\alpha<\kappa$. By
$M_\alpha\in M$, we have $[M_\alpha]^{<\kappa}\in M$. By 
$2^{<\kappa}=\kappa$, 
$[M_\alpha]^{<\kappa}$ has cardinality $\kappa$. Hence, by 
\Lemmaof{V-kappa-like},\,\assertof{a}, it follows that 
$[M_\alpha]^{\kappa}\subseteq M$. 
\smallskip\\ 
\assertof{2}\implies\assertof{1}: Suppose that $M\prec\calH_\chi$,
$\cardof{M}=\kappa$, $B\in M$ and $[M]^{<\kappa}\subseteq M$. Let
$M=\setof{m_\alpha}{\alpha<\kappa}$. Then we can construct inductively an 
increasing sequence $(M_\alpha)_{\alpha<\kappa}$ of elementary submodels of 
$M$ of cardinality less than $\kappa$ \st\ 
$M_\alpha,m_\alpha\in M_{\alpha+1}$ for every 
$\alpha<\kappa$. 
This is possible since at $\alpha$'th step of the inductive construction, we 
have that $M_\alpha$ is a subset of $M$ of cardinality less than $\kappa$. By
$[M]^{<\kappa}\subseteq M$, it follows that $M_\alpha\in M$. So by 
the downward L\"owenheim-Skolem Theorem, we can take 
$M_{\alpha+1}\prec M$ \st\ $\cardof{M_{\alpha+1}}<\kappa$ and 
$M_\alpha,\,m_\alpha\in M_{\alpha+1}$. 
At a limit $\gamma<\kappa$ we take $\bigcup_{\alpha<\gamma}M_\alpha$. 
Then $(M_\alpha)_{\alpha<\kappa}$ witnesses the 
$V_{\kappa}$-likeness of $M$.
\qedofLemma
\begin{Prop}\label{char-under-ch}
Assume that $\kappa$ is a regular cardinal \st\ $2^{<\kappa}=\kappa$. 
Then for a partially ordered structure $B$, 
the following are equivalent:\smallskip\\
\assert{1} $B$ has the $\kappa$-FN;\smallskip\\
\assert{2} For sufficiently large $\chi$ and for all $M\prec\calH_\chi$, if 
$B\in M$, $\cardof{M}=\kappa$ and $[M]^{<\kappa}\subseteq M$, then
$B\cap M\leq_\kappa B$ holds;\smallskip\\
\assert{3} Player II has a winning strategy in $\kappagame{\kappa}{B}$. 
\end{Prop}
\prf
\assertof{1}\equivto\assertof{2}: By \Lemmaof{Vomega-1-like} and 
\Propof{game}. \smallskip\\
\assertof{1}\implies\assertof{3} follows from \Propof{game}.\smallskip\\
\assertof{3}\implies\assertof{2}: Let $M$ be as in \assertof{2} and let 
$B\cap M=\setof{b_\alpha}{\alpha<\kappa}$. By \assertof{3}, there is a 
winning strategy $\tau\in M$ of Player II in $\kappagame{\kappa}{B}$. Let 
$(X_\alpha, Y_\alpha)_{\alpha<\kappa}$ be a play in 
$\kappagame{\kappa}{B}$ \st\ 
Player I chooses $X_\alpha$ so that $\cardof{X_\alpha}<\kappa$ and 
$b_\alpha\in X_\alpha$, and Player II played always according to 
$\tau$. Such a 
game is possible since, by $[M]^{<\kappa}\subseteq M$, at Player II's 
$\alpha$'th innings, she has $(X_0,Y_0,\ldots,X_\alpha)\in M$. Hence her 
move $Y_\alpha$ taken according to $\tau$ will be also an element of $M$. 
Since $\cardof{Y_\alpha}<\kappa$, 
$Y_\alpha$ is a subset of $M$. Now we have 
$\bigcup_{\alpha<\kappa}Y_\alpha=\bigcup_{\alpha<\kappa}X_\alpha=B\cap M$. 
Since $\tau$ was a winning strategy, we also have 
$B\cap M=\bigcup_{\alpha<\kappa}X_\alpha\leq_\kappa B$.\smallskip\\
\qedofProp
\indentafterqed
We can also consider the following variant of the game 
$\kappagame{\kappa}{B}$: for cardinals $\kappa$, $\lambda$ \st\ 
$\lambda\leq\kappa$ and 
a partially ordered structure $B$, $\kappalambdagame{\kappa}{\lambda}{B}$ is 
the game just like $\kappagame{\kappa}{B}$ except that 
Player II wins in $\kappalambdagame{\kappa}{\lambda}{B}$ if and only if
$\bigcup_{\alpha<\kappa}X_\alpha\leq_\lambda B$. 
As in Propositions \ref{characterization},\,\ref{game}, we can prove the 
implication 
\assertof{A}\implies\assertof{B}\implies\assertof{C}\implies\assertof{D} for 
the following assertions for regular $\kappa$, $\lambda$. \bigskip\\
\assert{A} For every sufficiently large $\chi$ and $M\prec\calH_\chi$ of 
cardinality $\kappa$ with $B\in M$, we have $B\cap M\leq_\lambda B$.\medskip\\
\assert{B} Player II has a simple winning strategy in 
$\kappalambdagame{\kappa}{\lambda}{B}$.\medskip\\ 
\assert{C} For every sufficiently large $\chi$ and $M\prec\calH_\chi$, if
$M$ is $V_{\kappa}$-like then $B\cap M\leq_\lambda B$.\medskip\\
\assert{D} Player II has a winning strategy in 
$\kappalambdagame{\kappa}{\lambda}{B}$.\medskip\\ 
\assert{E} For any sufficiently large $\chi$ and $M\prec\calH_\chi$ of 
cardinality $\kappa$ \st\ $B\in M$ and $[M^{<\kappa}\subseteq M$,
$B\cap M\leq_\lambda B$ holds.\bigskip\\
For the implication \assertof{A}\implies\assertof{B}, we fix an expansion
of $\calH_\chi$ by Skolem functions. For
$x\subseteq\calH_\chi$, let $\tilde{h}(x)$ be the Skolem hull of $x$.
We may take the Skolem hull operation so that $B\in\tilde{h}(\emptyset)$ 
holds. Player II then 
wins if she takes $Y_\alpha$ \st\ $X_\alpha\subseteq Y_\alpha$ and
$\tilde{h}(Y_\alpha)\cap B=Y_\alpha$ hold in each of her $\alpha$'th innings 
for $\alpha<\kappa$. 
By the same idea, we can also prove the equivalence of \assertof{B} 
and \assertof{C}, if we allow Player II to remember her last move 
in her simple winning strategy in \assertof{B}.
By \Lemmaof{Vomega-1-like}, we have 
\assertof{C}\equivto\assertof{D}\equivto\assertof{E} under 
$2^{<\kappa}=\kappa$. 
\section{Cardinal functions on \Bas\ with the weak Freese-Nation property}%
\label{card-functions}
In \cite{heindorf-shapiro} it is shown that, for any openly generated \Ba\ 
(i.e., \Ba\ with the FN), 
the cardinal functions (those studied in \cite{monk}, possibly except the 
topological density $d$) 
have the same value as for the free \Ba\ of the same cardinality, 
as follows obviously from \Lemmaof{aleph-1}. Later we shall see some more 
examples of \Bas\ with 
the WFN which behave quite differently from free \Bas\ \wrt\ cardinal 
functions. Nevertheless, there are some restrictions on the values of 
cardinal functions on \Bas\ with the WFN.
\begin{Prop}\label{omega-2}
For every partially ordered structure $B$, if $\kappa^++1$ or $(\kappa^++1)^*$ 
is (order isomorphic) embeddable into $B$ then $B$ does not have the 
$\kappa$-FN. In particular, for every \Ba\ with the $\kappa$-FN, we have 
${\rm Depth}(B)\leq\kappa$. 
\end{Prop}
\prf
Suppose that $\mapping{i}{\kappa^++1}{B}$ is an embedding (the case for 
$(\kappa^++1)^*$ can be handled similarly). Let $\mapping{j}{B}{\kappa^++1}$ be 
defined by 
\[ j(b)=\sup\setof{\alpha}{i(\alpha)\leq b}
\]\noindent
for $b\in B$. Then $j$ is order preserving and $j\circ i=id_{\kappa^++1}$ 
holds. 
Hence, by \Lemmaof{retract}, the following claim proves the proposition:
\begin{Claim}
$(\kappa^++1,\leq)$ does not have the $\kappa$-FN.
\end{Claim}
\prfofClaim
Player I wins a game in $\kappagame{\kappa}{\kappa^++1}$ if he chooses 
$X_\alpha$ at his $\alpha$'th move \st\ 
$\sup X_\alpha\setminus\smallsetof{\kappa^+}
	>\sup\bigcup_{\beta<\alpha}Y_\beta\setminus\smallsetof{\kappa^+}$ holds. 
By \Propof{game}, it follows that $(\kappa^++1,\leq)$ does not have 
the $\kappa$-FN. Note that for the implication 
\assertof{1}\implies\assertof{2} in \Propof{game} 
used here, we do not need the assumption of regularity of $\kappa$. 
\qedofClaim
\qedofProp 
\begin{Thm}\label{ind}
For a regular cardinal $\kappa$, if a \Ba\ $B$ has the 
$\kappa$-FN, $\lambda=\lambda^{<\kappa}$ and $X\subseteq B$ is of 
cardinality $>\lambda$ then there is an independent $Y\subseteq X$ of 
cardinality $>\lambda$. 
\end{Thm}
\prf
Essentially the same argument as the following one has been used in 
\cite[\S 4]{Sh:92}. 

Let $f$ be a $\kappa$-FN mapping on $B$. Let 
$(a_\delta)_{\delta<\lambda^+}$ be a sequence of elements of $X$ \st, 
letting $B_\delta$ be the closure of $\setof{a_\gamma}{\gamma<\delta}$ \wrt\ 
$f$ and the Boolean operations, $a_\delta\not\in B_\delta$ holds for every 
$\delta<\lambda^+$. 
By \Lemmaof{sigma-embed}, \assertof{c}, we have $B_\delta\leq_\kappa B$ for 
every $\delta<\lambda^+$. 
Let
$S=\setof{\delta<\lambda^+}{\cf(\delta)\geq\kappa}$.
For each $\delta\in S$, let $I_\delta$ and $J_\delta$ be 
cofinal subsets of $B_\delta\restr a_\delta$ and $B_\delta\restr -a_\delta$ 
respectively, both of cardinality less than $\kappa$. Let
\[h(\delta)=\tuple{I_\delta, J_\delta}.
\]\noindent 
By Fodor's lemma and $\lambda=\lambda^{<\kappa}$, there is a stationary 
$T\subseteq S$ \st\ $h\restr T$ is 
constant, say $h(\delta)=\tuple{I,J}$ for all $\delta\in T$. 
Let 
\[
\delta^*=\min\setof{\delta<\lambda^+}{I,\,J\subseteq B_\delta}.
\] 
\Wolog\ we may assume that $\delta^*<\delta$ holds for every $\delta\in T$. 
Let 
\[ L=\setof{b\in B_{\delta^*}}{b\not\leq i+j\mbox{ for all }i\in I,\, j\in J}.
\]\noindent
Then we have\medskip\\
\assert{1} $1\in L$ (since, by $a_\delta\not\in B_{\delta^*}$ for any 
$\delta\in T$, $I\cup J$ generates a proper ideal of $B_{\delta^*}$). In 
particular we have $L\not=\emptyset$;\\
\assert{2} If $b\in L$ and $k\in I\cup J$ then $b\cdot-k\in L$. \medskip\\
Now, by \assertof{1} above, the following claim shows that 
$\setof{a_\delta}{\delta\in T}$ is independent. Since 
$\cardof{T}=\lambda^+$, this proves the theorem.
\begin{Claim}
If $b\in L$ and $p$ is an elementary product over 
$a_{\delta_0}$,\ldots,$a_{\delta_{n-1}}$ for $\delta_i\in T$ \st\ 
$\delta_0<\cdots<\delta_{n-1}$ 
(i.e.\ $p$ is of the form 
$(a_{\delta_0})^{\tau_0}\cdot\,\cdots\,\cdot(a_{\delta_{n-1}})^{\tau_{n-1}}$ 
for some $\tau_i\in 2$, $i<n$) then $b\cdot p\not= 0$. 
Here, for a \Ba\ 
$B$, $b\in B$ and $i\in 2$, we define $(b)^i$ by:
\[ (b)^i=
\left\{\,
\begin{array}{@{}l@{}l}
b\quad&\mbox{; if }i=1,\\
-b\quad&\mbox{; if }i=0.
\end{array}
\right.
\]\noindent
\end{Claim}
\prfofClaim
By induction on $n$. For $n=0$, this is trivial since $0\not\in L$. Assume 
that the claim holds for $n$. Let $\delta_0$,\ldots,$\delta_n\in T$ be \st\ 
$\delta_0<\cdots\delta_{n-1}<\delta_n$ and let $p$ be an arbitrary 
elementary product over $a_{\delta_0}$,\ldots,$a_{\delta_{n-1}}$. 
Let $b\in L$. By the 
induction hypothesis, we have $b\cdot p\not=0$. We have to show that 
$b\cdot p\cdot a_{\delta_n}\not=0$ and $b\cdot p\cdot-a_{\delta_n}\not=0$. 
Toward a contradiction, assume that $b\cdot p\cdot a_{\delta_n}=0$ holds. Then 
$b\cdot p\leq -a_{\delta_n}$. Since $b\cdot p\in B_{\delta_n}$, we can find 
$j\in J$ \st\ $b\cdot p\leq j$. Hence $(b\cdot-j)\cdot p=0$. Since 
$b\cdot-j\in L$ by \assertof{2} above, this is a contradiction to the 
induction hypothesis. Similarly, from $b\cdot p\cdot-a_{\delta_n}=0$, it 
follows that $(b\cdot-i)\cdot p=0$ for some $i\in I$ which again is a 
contradiction to \assertof{2}.
\qedofClaim\\
\qedofThm
\indentafterqed
The next corollary gives a positive answer to a problem 
by L.\ Heindorf. 
\begin{Cor}\label{ind-ineq}
For a regular cardinal $\kappa$, if a \Ba\ 
$B$ has the $\kappa$-FN then $\cardof{B}\leq{\rm Ind}(B)^{<\kappa}$. In 
particular, 
for an openly generated \Ba\ $B$, $\cardof{B}= {\rm Ind}(B)$ 
holds (\cite{heindorf-shapiro}). 
For 
a \Ba\ $B$ with the WFN, we have 
$\cardof{B}\leq{\rm Ind}(B)^{\aleph_0}$. 
\end{Cor}
\prf
Assume that $B$ has the $\kappa$-FN but 
$\cardof{B}>{\rm Ind}(B)^{<\kappa}$ holds. Then by 
${\rm Ind}(B)^{<\kappa}=({\rm Ind}(B)^{<\kappa})^{<\kappa}$, we have 
${\rm Ind}(B)>{\rm Ind}(B)^{<\kappa}$ by \Thmof{ind}. But this is impossible. 
\qedofCor
\noindentafterqed
In \Corabove\ the equality is attained for every regular $\kappa$. For the 
case of $\kappa=\aleph_1$, 
the simplest example to see this would be $Intalg(\reals)$ 
(see \Propof{intalg} below). 

The following corollary is an immediate consequence of \Thmof{ind}:
\begin{Cor}{\rm(\cite{heindorf-shapiro} for $\kappa\leq\aleph_1$)}\label{Lutz}
Let $\kappa$ be a regular cardinal. If a \Ba\ $B$ has the $\kappa$-FN, then 
$c(B)\leq2^{<\kappa}$ and ${\rm Length}(B)\leq 2^{<\kappa}$. \qed 
\end{Cor}
\section{Interval algebras and power set algebras } \label{Intalg}
\begin{Prop}\label{intalg}
\mbox{}\\ 
\assert{a} For $\rho\in\Ord$, $Intalg(\rho)$ has 
the $\kappa$-FN if and only if $\rho<\kappa^+$.\smallskip\\
\assert{b}$Intalg(\reals)$ has the WFN.\smallskip\\
\assert{c} For a totally ordered set $X$, $Intalg(X)$ has the $\kappa$-FN 
if and only if $X$ has the $\kappa$-FN.\smallskip\\ 
\assert{d} Assume that $\kappa$ is a regular cardinal. For a linearly 
ordered set $X$, if $Intalg(X)$ (hence, 
by \assertof{c}, also $X$) has the $\kappa$-FN then 
$\cardof{X}\leq 2^{<\kappa}$. 
\end{Prop}
\prf
For $b\in Intalg(X)$, let $ep(b)$ be the set of end points of 
$b$, i.e.
\[ 
ep(b)=\setof{x_i}{i<2n}. 
\]\noindent
where $b=\dot{\bigcup}_{j<n}[x_{2j},x_{2j+1})$ in the standard 
representation. \medskip\\
\assertof{a}: For $\rho<\kappa^+$, $Intalg(\rho)$ has cardinality less or 
equal to $\kappa$. Hence, by \Lemmaof{aleph-1}, $Intalg(\rho)$ 
has the $\kappa$-FN. 
If $\rho\geq\kappa^+$, $Intalg(\rho)$ does not have the $\kappa$-FN 
by \Propof{omega-2}.\smallskip\\ 
\assertof{b}: 
For all $b\in Intalg(\reals)$, the mapping defined by
\[g(b)=\setof{c\in Intalg(\reals)}{ep(c)\subseteq\rationals
\cup ep(b)}.\]\noindent 
is a WFN-mapping on $Intalg(\reals)$.\smallskip\\
\assertof{c}: 
If 
$\mapping{f}{X}{[X]^{<\kappa}}$ is a $\kappa$-FN mapping on $X$ then 
$\mapping{g}{Intalg(X)}{[Intalg(X)]^{<\kappa}}$ defined by 
\[ g(b)=\setof{c\in Intalg(X)}{ep(c)\subseteq\bigcup f[ep(b)]}
\]\noindent
is a $\kappa$-FN mapping on $Intalg(X)$. Conversely, if $g$ is a 
$\kappa$-FN mapping on 
$Intalg(X)$, then $\mapping{f}{X}{[X]^{<\kappa}}$ defined by
\[ f(x)=\bigcup\setof{ep(b)}{b\in g((-\infty,\,x))}
\]\noindent
is a $\kappa$-FN mapping on $X$.\medskip\\
\assertof{d}: We have ${\rm Ind}(Intalg(X))=\aleph_0$ (see e.g.\ Corollary 
15.15 in \cite{koppelbook}) Hence, by \Corof{ind-ineq}, 
$\cardof{X}\leq\cardof{Intalg(X)}\leq 2^{<\kappa}$ holds if $X$ (or 
equivalently $Intalg(X)$) has the $\kappa$-FN.
\qedofProp
\noindentafterqed
Note that \Lemmaabove,\,\assertof{c} 
is not true for tree algebras: e.g., for any cardinal 
$\lambda$ \st\ $\lambda>2^{\aleph_0}$, the tree $(\kappa,\emptyset)$ 
has the WFN but $Treealg((\kappa,\emptyset))$ does not by \Corof{Lutz}.
\par
For any sets $x$, $y$ we say that $x$ is a subset of $y$ modulo 
$<\kappa$ (notation: $x\subseteq_{<\kappa}y$) if 
$\cardof{x\setminus y}<\kappa$ holds. The following lemma is 
well-known: 
\begin{Lemma}\label{tower-in-omega-2}Suppose that $\kappa$ is a regular 
uncountable 
cardinal.\smallskip\\ 
\assert{a} If $(u_\alpha)_{\alpha<\kappa}$ is a sequence of non-stationary 
subsets of $\kappa$ then there exists a non-stationary 
$u\subseteq\kappa$ \st\ $u_\alpha\subseteq_{<\kappa}u$ holds for all 
$\alpha<\kappa$. \smallskip\\
\assert{b} There exists a strictly $\subseteq_{<\kappa}$-increasing sequence 
of elements of $\powersetof{\kappa}$ of order type $\kappa^+$.
\end{Lemma}
\prf
\assertof{a}: For each $\alpha<\kappa$, let 
$c_\alpha\in\powersetof{\kappa}$ be a club subset of $\kappa$ \st\ 
$u_\alpha\cap c_\alpha=\emptyset$. 
Let $(\beta_\alpha)_{\alpha<\kappa}$ be a strictly, continuously 
increasing sequence in 
$\kappa$ \st\ $\beta_\alpha\in\bigcap_{\delta<\alpha}c_\delta$ for every 
$\alpha<\kappa$. Then 
$\kappa\setminus\setof{\beta_\alpha}{\alpha<\kappa}$ is as 
desired.\smallskip\\
\assertof{b}: We can construct a $\subseteq_{<\kappa}$-increasing sequence 
$(u_\alpha)_{\alpha<\kappa}$ of elements of $\powersetof{\kappa}$ 
inductively so 
that $u_\alpha$ is non-stationary for every $\alpha<\kappa$:
for a successor step let $u_{\alpha+1}$ be the union of $u_\alpha$ and any 
non-stationary subset of $\kappa\setminus u_\alpha$ of cardinality less than 
$\kappa$. For a limit $\delta<\kappa^+$ with 
$\cf(\delta)=\lambda<\kappa$, 
we choose increasing 
$(\delta_\beta)_{\beta<\lambda}$ \st\ 
$\delta=\bigcup_{\beta<\lambda}\delta_\beta$, and let 
$u_\delta=\bigcup_{\beta<\lambda}u_{\delta_\beta}$. For limit 
$\delta<\kappa^+$ with $\cf(\delta)=\kappa$, we can take an appropriate 
$u_\delta$ using \assertof{a}. 
\qedofLemma
\begin{Prop}\label{P(omega-n)}Suppose that $\kappa$ is a regular uncountable 
cardinal. Then\smallskip\\
\assert{a}$\powersetof{\kappa}/[\kappa]^{<\kappa}$ does not have the 
$\kappa$-FN. \smallskip\\
\assert{b} $\powersetof{\kappa}$ does not have the $\kappa$-FN.
\end{Prop}
\prf 
\assertof{a}: By \Lemmaof{tower-in-omega-2},\,\assertof{b}, 
$(\kappa^+,\,\leq)$ is 
embeddable into $\powersetof{\kappa}/[\kappa]^{<\kappa}$. Hence, 
by \Propof{omega-2}, $\powersetof{\kappa}/[\kappa]^{\kappa}$ does not have the 
$\kappa$-FN. \smallskip\\
\assertof{b}: 
Let $\chi$ be sufficiently large and let $M\prec\calH_\chi$ be 
$V_{\kappa}$-like \st\ $\kappa\in M$ (and hence also 
$\powersetof{\kappa}\in M$). Let
$(M_\alpha)_{\alpha<\kappa}$ be an increasing 
sequence of elementary submodels of $M$ of cardinality less than $\kappa$ 
witnessing the 
$V_{\kappa}$-likeness of $M$. We construct a 
sequence 
$(u_\alpha)_{\alpha<\kappa}$ of non-stationary subsets of $\kappa$ 
inductively \st\ $u_\alpha\in M_{\alpha+1}$, 
$u\subseteq_{<\kappa}u_\alpha$ 
and $\cardof{u_\alpha\setminus u}=\kappa$ hold for every non-stationary 
$u\in\powersetof{\kappa}\cap M_\alpha$. 
This is possible since $M_\alpha\in M_{\alpha+1}$. By 
\Lemmaof{tower-in-omega-2},\,\assertof{a}, there is a non-stationary 
$u^*\in\powersetof{\kappa}$ \st\ $u_\alpha\subseteq_{<\kappa}u^*$ 
holds for every $\alpha<\kappa$. Clearly 
$(\powersetof{\kappa}\cap M)\restr u^*$ is not generated by any subset of 
cardinality less than $\kappa$. 
Hence, by \Propof{game}, it follows that $\powersetof{\kappa}$ does not 
have the $\kappa$-FN. 
\qedofProp
\indentafterqed
By \Propof{P(omega-n)}, it follows that $\powersetof{\omega_1}$ does not 
have the WFN. In contrast to this, the 
statement ``$\powersetof{\omega}$ has the WFN'' is independent 
from \ZFC\ or even from \ZFC\ $+$ $\neg$\CH. 
For $x\in\powersetof{\omega}$, let us denote by $[x]$ the equivalence 
class of $x$ modulo fin (the ideal of finite subsets of $\omega$).
\begin{Lemma}\label{fin}
$\powersetof{\omega}$ has the WFN if and only if 
$\powersetof{\omega}/fin$ has the WFN.
\end{Lemma}
\prf
If $g$ is a WFN mapping on 
$\powersetof{\omega}$, then 
$\mapping{g'}{\powersetof{\omega}/fin}{%
[\powersetof{\omega}/fin]^{\leq\aleph_0}}$ 
defined by
\[ g'([x])
=\bigcup\setof{\setof{[y]}{y\in g(z)}}{z\in\powersetof{\omega},\,
[z]=[x]}
\]\noindent
for $x\in\powersetof{\omega}$, is a WFN mapping on $\powersetof{\omega}/fin$.

If $f$ is a WFN mapping on $\powersetof{\omega}/fin$, then 
$\mapping{f'}{\powersetof{\omega}}{[\powersetof{\omega}]^{\leq\aleph_0}}$ 
defined by
\[ f'(x)=\setof{z}{[z]\in f([x])}
\]\noindent
for $x\in\powersetof{\omega}$, is a WFN mapping on $\powersetof{\omega}$.
\qedofLemma
\begin{Prop}\label{P(omega)}
\assert{a}{\rm(\CH)} $\powersetof{\omega}$ has the WFN.\smallskip\\
\assert{b} In the generic extension of a model of\/ \ZFC $+$ \CH\ by adding 
less than $\aleph_\omega$ many 
Cohen reals (by standard Cohen forcing), 
$\powersetof{\omega}$ still has the WFN. In 
particular the assertion 
``\/$\MA(\mbox{\it Cohen})$ $+$ $\neg$\CH\ $+$ 
$\powersetof{\omega}$ has the WFN'' is consistent. 
Here $\MA(\mbox{\it Cohen})$ stands for Martin's axiom restricted to 
\pos\ of the form $\Fn(\kappa,2)$.
\smallskip\\
\assert{c} If ${\bf b}\geq\aleph_2$, then $\powersetof{\omega}$ does not 
have the WFN. 
\end{Prop}
\prf
\assertof{a}: Since $\cardof{\powersetof{\omega}}=\aleph_1$ under \CH, 
the claim follows from \Lemmaof{aleph-1}.\smallskip\\
\assertof{b}: 
This follows from \Thmof{cohen-model} below.
\smallskip\\
\assertof{c}: By \Lemmaof{fin}, it is enough to show that 
$\powersetof{\omega}/fin$ does not have the WFN. But under 
${\bf b}\geq\aleph_2$, $(\omega_2,\leq)$ can be embedded into 
$\powersetof{\omega}/fin$. Hence, by \Propof{omega-2}, 
$\powersetof{\omega}/fin$ does not have the WFN.
\qedofProp
\noindentafterqed
Note that, by \Propabove,\,\assertof{c}, the 
statement ``$\powersetof{\omega}$ has the WFN'' is not consistent 
with \MA({\it$\sigma$-centered\/}) $+$ $\neg$\CH. 
\section{Complete \Bas}
\begin{Lemma}\label{P(kappa)embedded}
For \Bas\ $A$, $B$ and regular $\kappa$, if $A\leq B$, $A$ is 
complete (but not necessarily a complete subalgebra of $B$) and $B$ has the 
$\kappa$-FN, then $A$ also has the $\kappa$-FN. 
\end{Lemma}
\prf
By Sikorski's Extension Theorem, there is a homomorphism $j$ from $B$ to 
$A$ \st\ $j\restr A=id_A$. Hence the claim of the lemma follows from 
\Lemmaof{retract}. 
\qedofLemma
\begin{Thm}\label{complete-Ba}
\assert{a} If ${\bf b}\geq\aleph_2$, then no infinite complete \Ba\ has the 
WFN.\smallskip\\ 
\assert{b} 
For regular $\kappa$, no \cBa\ without the $\kappa$-cc has the 
$\kappa$-FN.\smallskip\\ 
\assert{c} If $\kappa$ is regular and $2^{<\kappa}=\kappa$, then every 
$\kappa$-cc \cBa\ has the $\kappa$-FN.
\end{Thm}
\prf
\assertof{a}: If $B$ is complete, then $\powersetof{\omega}$ is embeddable 
into 
$B$. Hence, by \Propof{P(omega)},\,\assertof{c} and \Lemmaof{P(kappa)embedded}, 
$B$ does not have the WFN.\smallskip\\
\assertof{b}: If $B$ is complete and not $\kappa$-cc, then 
$\powersetof{\kappa}$ is embeddable into $B$. Hence, by 
\Propof{P(omega-n)},\,\assertof{b} and \Lemmaof{P(kappa)embedded}, 
$B$ does not have the $\kappa$-FN.\smallskip\\
\assertof{c}: Assume that $B$ has the $\kappa$-cc. 
Let $\chi$ be sufficiently large and let $M\prec\calH_\chi$ be \st\ 
$\cardof{M}=\kappa$, $B\in M$ and $[M]^{<\kappa}\subseteq M$ hold. 
Then $B\cap M$ is a $\kappa$-complete subalgebra of $B$. By the $\kappa$-cc 
of $B$ it follows that $B\cap M$ is complete subalgebra of $B$. Hence 
$B\cap M\leq_\rc B$ holds. 
By \Propof{char-under-ch} it follows that $B$ has the $\kappa$-FN.
\qedofThm
\indentafterqed
By \assertof{b} and \assertof{c} in \Thmabove, under \CH,
a complete \Ba\ $B$ has the WFN if and only if $B$ satisfies the 
\ccc. 
This assertion still holds partially in the model obtained by adding Cohen 
reals to a model of \CH:
\begin{Thm}\label{cohen-model}
Suppose that $V\models{\ZFC+\CH}$. Let $\lambda$ be a cardinal in 
$V$ \st\ $V\models{\lambda<\aleph_\omega}$. Let $P=\Fn(\lambda, 2)$ and let  
$G$ be a $V$-generic filter $G$ over $P$. For $B\in V[G]$ \st\ 
$V[G]\models{B\xmbox{ \it is \ccc\ \cBa}}$, if either\smallskip\\
\assert{a} $V[G]\models{\cardof{B}<\aleph_\omega}$ or \\
\assert{b} there is a \Ba\ $A\in V$ \st\ 
$B=\overline{A}^{V[G]}$\smallskip\\
$($where $\overline{A}^{V[G]}$ denotes the completion 
of $A$ in $V[G]$\relax$)$, then we have:
\[ V[G]\models{B \mbox{ has the WFN}}.
\]\noindent
\end{Thm}
For the proof of \Thmabove\ we need the following lemma. 
\begin{Lemma}\label{Hilfslemma}
Let $V$ be a ground model and $P=\Fn(S,2)$ for some $S\in V$. Let $G$ be 
$V$-generic over $P$ and $A\in V$ be \st\ 
\[ V\models{A\mbox{ is a \ccc\ \cBa}}.
\]\noindent
Then we have:\\
\assert{a} if $\calS\in V$ is \st\ 
\[ 
V\models{\calS\mbox{ is a }\sigma\mbox{-directed family of subsets of }S
	\mbox{ and }\bigcup\calS=S}
\]\noindent
then $\overline{A}^{V[G]}=\bigcup\setof{\overline{A}^{V[G_T]}}{T\in\calS}$ 
where $G_T=G\cap \Fn(T,2)$;\smallskip\\
\assert{b} For a \Ba\ $B$ in $V[G], $ if $V[G]\models{A\leq B}$, then 
$V[G]\models{ A\leq_\sigma B}$.
\end{Lemma}
\prf
\assertof{a}: If $b\in\overline{A}^{V[G]}$ then, by the \ccc\ of $A$, there 
is a countable $X\subseteq A$ ($X\in V[G]$) %
\st\ $b=\Sigma^{\overline{A}^{V[G]}}X$. 
Hence, by the \ccc\ of $\Fn(S,2)$ (in $V$), there is a name 
$\dotb$ of $b$ in which only countably many elements of $\Fn(S,2)$ 
appear. Let $T\in\calS$ be \st\ every element of  
$\Fn(S,2)$ appearing in $\dotb$ is contained in $\Fn(T,2)$. 
Then $b\in V[G_T]$, hence $b\in \overline{A}^{V[G_T]}$.\smallskip\\
\assertof{b}: 
It is enough to show that $V[G]\models{A\leq_\sigma\overline{A}^{V[G]}}$. 
Let $b\in\overline{A}^{V[G]}$. Then, as in the proof 
of \assertof{a}, there exists a countable $X\subseteq S$ (in $V$) %
\st\ $b\in V[G_X]$ where $G_X=G\cap\Fn(X,2)$. Let $\dotb$ be an 
$\Fn(X,2)$ name of $b$. In $V[G]$, let 
\[ I=\setof{\Sigma^A\setof{c\in A}{p\forces{\Fn(X,2)}{c\leq\dotb}}}{p\in G_X}.
\]\noindent
Then $I$ is countable and generates $A\restr b$.
\qedofLemma
\noindentafterqed
\prfof{\Thmof{cohen-model}} 
First, let us consider the case \assertof{b}. 
It is enough to show the following assertion for all \ccc\ \Bas\ $A$ in $V$  
and for all $n\in\omega$: 
\begin{assertion}{$(*)_n$}\it
If $H$ is $V$-generic over $\Fn(\aleph_n,2)$ then $\overline{A}^{V[H]}$ has 
the WFN in $V[H]$.
\end{assertion}
For $n\leq 1$, $V[H]$ still satisfies the \CH. Hence, by 
\Thmof{complete-Ba},\,\assertof{c}, $\overline{A}^{V[H]}$ has the WFN in 
$V[H]$. Now suppose that $n>1$ and we have shown $(*)_m$ for all $m<n$. Let 
$H$ be a $V$-generic filter over $\Fn(\aleph_n,2)$. For each 
$\alpha<\aleph_n$, let $H_\alpha=H\cap\Fn(\alpha,2)$. Then by the induction 
hypothesis and \Lemmaof{Hilfslemma},\,\assertof{a} and \assertof{b}, the 
sequence $(\overline{A}^{V[H_\alpha]})_{\alpha<\aleph_n}$ satisfies the 
conditions of \Lemmaof{non-conti-chain}. Hence 
$\overline{A}^{V[H]}=\bigcup_{\alpha<\aleph_n}\overline{A}^{V[H_\alpha]}$ has 
the WFN in $V[H]$. \medskip\par
Now suppose that $B$ is a \ccc\ \cBa\ in $V[G]$ \st\ 
$V[G]\models{\cardof{B}<\aleph_\omega}$. 
\Wolog\ we may assume that the underlying set of $B$ is a cardinal 
$\kappa<\aleph_\omega$. Let $\dot\leq$ be a 
$P$-name of the partial ordering of $B$. In $V$, let 
$u_{\alpha,\beta}\subseteq\setof{p\in P}{%
p\xmbox{ decides }\alpha\mathrel{\dot\leq}\beta}$ be maximally pairwise 
disjoint for each $\alpha,\beta\in\kappa$ (more precisely the family 
$(u_{\alpha,\beta})_{\alpha,\beta\in\kappa}$ should be taken in $V$).
Since $P$ satisfies the \ccc, we have 
$\cardof{u_{\alpha,\beta}}\leq\aleph_0$. Further in $V$, let 
\[ 
\begin{array}{@{}l@{}l}
I=\setof{(r,s)}{%
	& r\in[\kappa]^{\aleph_1},\,s\in[\lambda]^{\aleph_1},\\ 
	& \forces{P}{r\mbox{ is closed \wrt\ Boolean operations}}\\
	& \mbox{ and }u_{\alpha,\beta}\subseteq\Fn(s,2)
		\mbox{ for all }\alpha,\beta\in r
}.
\end{array}
\]\noindent
For $(r,s),\,(r',s')\in I$, let 
\[ (r,s)\leq(r',s')\ \ \Leftrightarrow\ \ r\subseteq r'\mbox{ and }
s\subseteq s'.
\]\noindent
By the definition of $I$, we have $B\cap r\in V[G_s]$ for $(r,s)\in I$ where 
$B\cap r$ is the subalgebra of $B$ with the underlying set $r$ and, 
as before, $G_s=G\cap\Fn(s,2)$. Let 
$B_{(r,s)}=\overline{B\cap r}^{V[G_s]}$ and let $\dotB_{(r,s)}$ be its 
$\Fn(s,2)$-name. 
Since $V[G_s]\Models{\CH}$, we have 
$V[G_s]\models{\cardof{B_{(r,s)}}=\aleph_1}$ for every $(r,s)\in I$. 

The rest of the proof is modeled after the proof of 
\assertof{4}\implies\assertof{1} of \Propof{characterization}. 
For directed $I'\subseteq I$, let $B_{I'}=\bigcup_{(r,s)\in I'}B_{(r,s)}$. 
\addtocounter{Thm}{-1}
\begin{Claim}
If $I'\subseteq I$ is \st\ $I'\in V$ and 
$V\models{I'\xmbox{ is }\sigma\xmbox{-directed}}$, then 
$V[G]\models{B_{I'}\leq_\sigma B}$. 
\end{Claim}
\prfofClaim
Otherwise there is some $b\in B$ \st\ $B_{I'}\restr b$ is not countably 
generated in $V[G]$. Let $\dotb$ be a $P$-name of $b$. In $V$, we can 
construct an increasing sequence $(r_\alpha,s_\alpha)_{\alpha<\omega_1}$ of 
elements of $I'$ \st\ 
\[ \forces{P}{\dotB_{(r_\alpha,s_\alpha)}\restr\dotb\mbox{ does not generate }
\dotB_{(r_{\alpha+1},s_{\alpha+1})}\restr\dotb}. 
\]\noindent
Note that this is possible because of the \ccc\ of $P$. Let 
$r^*=\bigcup_{\alpha<\omega_1}r_\alpha$ and 
$s^*=\bigcup_{\alpha<\omega_1}r_\alpha$. Then $(r^*,s^*)\in I$. By an 
argument similar to the proof of \Lemmaof{Hilfslemma},\,\assertof{a}, we can 
show that $B_{(r^*,s^*)}=\bigcup_{\alpha<\omega_1}B_{(r_\alpha,s_\alpha)}$. 
Hence, by the construction of $(r_\alpha,s_\alpha)$, $\alpha<\omega_1$, the 
ideal $B_{(r^*,s^*)}\restr b$ is not countably generated. This is a 
contradiction to \Lemmaof{Hilfslemma},\,\assertof{b}.
\qedofClaim\medskip\\
For $X\in\powersetof{\kappa}\cap V$ and $Y\in\powersetof{\lambda}\cap V$, 
let 
\[ I_{X,Y}=\setof{(r,s)\in I}{r\subseteq X,\,s\subseteq Y}.
\]\noindent
Note that $I_{X,Y}$ is $\sigma$-directed. By induction on 
$\mu\leq\max(\kappa,\lambda)$, we can show easily that
\begin{Claim}
For $X\in\powersetof{\kappa}\cap V$ and $Y\in\powersetof{\lambda}\cap V$, if 
$\cardof{X}$, $\cardof{Y}\leq\max(\kappa,\lambda)$ then $B_{I_{X,Y}}$ has 
the WFN in $V[G]$. \qed
\end{Claim}
Since $I=I_{\kappa,\lambda}$, it follows from 
\Lemmaof{Hilfslemma},\,\assertof{a} that $B_{I_{\kappa,\lambda}}=B_I=B$. 
Hence, by \Claimabove, $B$ has the WFN in $V[G]$.
\qedof{\Thmof{cohen-model}}
\addtocounter{Thm}{1}
\begin{Cor}
The assertion ``every Cohen algebra has the WFN'' is consistent with 
{\rm\MA({\it Cohen})} $+$ $\neg\CH$. \qed
\end{Cor}
\section{$L_{\infty\kappa}$-free \Bas}\label{Linftykappa-free}
A Boolean algebra $B$ is called {\em$L_{\infty\kappa}$-free}\/ if 
$B\equiv_{L_{\infty\kappa}}\Fr\kappa$, i.e.\ if 
$B$ is elementary equivalent to $\Fr\kappa$ in the infinitary logic 
$L_{\infty\kappa}$. $B$ is 
{\em$L_{\infty\kappa}$-projective}\/ if $B\equiv_{L_{\infty\kappa}}C$ for some 
projective $C$. It is easily seen that if $B$ is 
$L_{\infty\kappa}$-projective then $B\oplus\Fr\kappa$ is 
$L_{\infty\kappa}$-free. 
In \cite{FuKoTa} it is shown 
that for every $\kappa$, there exists an $L_{\infty\aleph_1}$-free \Ba\ $B$ 
which does not satisfy the $\kappa$-cc. By \Corof{Lutz}, it follows that
\begin{Prop}
For any $\kappa$, there exists an $L_{\infty\aleph_1}$-free \Ba\ 
which does not have the $\kappa$-FN. \qed
\end{Prop}
On the other hand, we show in \Corof{L-infty-aleph-1} below that 
every ccc $L_{\infty\aleph_1}$-free \Ba\ has the WFN. 
Let us begin with the 
following lemma:

\begin{Lemma}\label{union-of-rc}
If $B$ satisfies the \ccc\ and $(A_\alpha)_{\alpha<\kappa}$ is an 
increasing sequence of relatively complete subalgebras of $B$ with 
$\cf(\kappa)>\omega$ then $\bigcup_{\alpha<\kappa}A_\alpha\leq_\rc B$ holds. 
\end{Lemma}
\prf
Suppose not. Then there is some $b\in B$ without projection on 
$\bigcup_{\alpha<\kappa}A_\alpha$. Then, for 
$\mu=\cf(\kappa)$, we can construct an increasing sequence of 
ordinals 
$(\gamma_\beta)_{\beta<\mu}$ \st\ $\gamma_\beta<\kappa$ and 
$p^B_{A_{\gamma_\beta}}(b)<p^B_{A_{\gamma_{\beta+1}}}(b)$ holds for 
every $\beta<\mu$. But this is impossible since $B$ satisfies the \ccc.
\qedofLemma
\begin{Prop}\label{ccc-rc}
If a \Ba\ $B$ satisfies the \ccc\ and 
$\setof{C}{C\leq_\rc B,\,\cardof{C}=\aleph_0}$ is cofinal in 
$[B]^{\aleph_0}$ then $B$ has the WFN.
\end{Prop}
\prf
By \Propof{game}, it is enough to show that Player II has a simple winning 
strategy in $\kappagame{\omega_1}{B}$. By \Lemmaof{union-of-rc}, Player II 
wins if he 
takes $Y_\alpha\geq X_\alpha$ in his $\alpha$'th move \st\ 
$Y_\alpha\leq_\rc B$ holds (actually this is a simple winning strategy of 
Player II in $\kappalambdagame{\omega_1}{\aleph_0}{B}$).
\qedofProp
\begin{Cor}\label{L-infty-aleph-1}
Every \ccc\ $L_{\infty\aleph_1}$-projective \Ba\ $B$ has the WFN. 
\end{Cor}
\prf
Let $B$ be a \ccc\ $L_{\infty\aleph_1}$-projective \Ba. 
The statement ``$\setof{C}{C\leq_\rc B,\,\cardof{C}=\aleph_0}$ is cofinal in 
$[B]^{\aleph_0}$'' can be formulated in $L_{\infty\aleph_1}$ and is true in 
any projective \Ba. Hence it is also true in $B$. 
By \Propof{ccc-rc}, 
it follows that $B$ has the WFN. 
\qedofCor
\indentafterqed
Since the \ccc\ is expressible in $L_{\infty\aleph_2}$ and it is true in any 
projective \Ba, every $L_{\infty\aleph_2}$-projective \Ba\ satisfies the ccc. 
Hence, by \Corabove, every $L_{\infty\aleph_2}$-projective \Ba\ has the WFN. 
Under Axiom R we can obtain a stronger result. 
Recall that Axiom R is the following statement:
\begin{assertion}{(Axiom R):}\it
For any $\lambda\geq\aleph_2$ with uncountable cofinality, stationary 
$\calS\subseteq[\lambda]^{\aleph_0}$ and a 
$\calT\subseteq[\lambda]^{\aleph_1}$ which is closed under union of 
increasing chains of order type $\omega_1$, there exists an 
$X\in\calT$ \st\ $\calS\cap[X]^{\aleph_0}$ is stationary in $[X]^{\aleph_0}$.
\end{assertion}
Axiom R follows from \MM\ (\cite{beaudoin}) but its consistency 
with \CH\ can be also shown under a supercompact cardinal.  
The following theorem was proved in \cite{fuchino2} (see also 
\cite{appendix}): 
\begin{Thm}
{\rm(Axiom R)} Every $L_{\infty\aleph_2}$-projective \Ba\ $B$ has the FN. \qed
\end{Thm}
The theorem above is not provable in \ZFC:  
under $V=L$, we can obtain a counter-example to \Thmabove\ (\cite{appendix}). 
In contrast to \Corof{L-infty-aleph-1}, we have the following:
\begin{Prop}
For any cardinal $\kappa$, there is a subalgebra of $\Fr\kappa^+$ without the 
$\kappa$-FN.
\end{Prop}
\prf
The topological dual to the \Ba\ $B$ below is considered in Engelking
\cite{engelking} to show that there exists a non-projective subalgebra of a
free \Ba\ (in the language of topology, this means that there exists a 
dyadic space which is not a Dugundji space). 

Let $X$ be a set of cardinality $\kappa^+$. We shall show that there is a
subalgebra of $\Fr X$ without the WFN. 
Let $U_1$ and $U_2$ be the ultrafilters of $\Fr X$ generated by $X$ and 
$\setof{-x}{x\in X}$ respectively. 
Let 
\[ B=\setof{b\in\Fr X}{b\in U_1\,\Leftrightarrow\, b\in U_2}.
\]\noindent
Clearly $B$ is a subalgebra of $\Fr X$. We claim that $B$ does not have the
$\kappa$-FN. 
For $Y\subseteq X$, let $B_Y=B\cap\Fr Y$. 
\begin{Claim}
For every $Y\in[X]^{\kappa}$, $B_Y$
is not a $\kappa$-subalgebra of $B$. 
\end{Claim}
\prfofClaim
Let $x_0\in Y$ and let $x_1$, $x_2$ be two distinct elements of $X\setminus Y$.
Let 
\[ b=x_0+x_1+-x_2. 
\]\noindent
Since
$b\in U_1$ and $b\in U_2$, we have $b\in B$. Let $I=B_Y\restr b$. We show
that $I$ is not $<\kappa$-generated: let $J$ be any subset of
$I$ of cardinality less than $\kappa$. For $c\in J$, we have $c\leq x_0$. 
Since $x_0$ is not an element of 
$B$, $c$ is strictly less than $x_0$. Let $Y'$ be a subset of
$Y$ of cardinality less than $\kappa$ \st\ $J\subseteq B_{Y'}$. Let $y_1$, 
$y_2$ be 
two distinct elements of 
$Y\setminus Y'$. Let
\[ d=x_0\cdot y_1\cdot-y_2. 
\]\noindent
Then we 
have $d\not\in U_1$, $d\not\in U_2$ and $d\leq x_0$. Hence $d\in B_Y\restr b$.
But $d$ is incomparable with every non-zero element of $J$.
\qedofClaim
\noindentafterqed
Since there are club many $C\in[B]^{\kappa}$ of the form $C=B_Y$ for 
$Y\in[X]^{\kappa}$, it follows that $B$ does not have the $\kappa$-FN by 
\Propof{characterization}. Note that the implication 
\assertof{1}\implies\assertof{3} in 
\Propof{characterization} used here does not require the assumption of 
regularity of $\kappa$.
\qedofProp
\indentafterqed
A subalgebra of an openly generated \Ba\ $B$ is openly generated 
(i.e.\ has the FN) if and only if $B$ has the Bockstein separation 
property (\cite{heindorf-shapiro}). 
\begin{Problem}
Is there a \Ba\ with the Bockstein separation property but without the WFN? 
\end{Problem}
In \cite{appendix}, the following partial answer to the problem is given:
{\em if $B$ is stable and satisfies the ccc and the Bockstein separation 
property, 
then $B$ has the WFN.} Here a \Ba\ $B$ is said to be stable if, for 
any countable subset $X$ of $B$, only countably many types over $X$ are 
realized in $B$. 


\end{document}